\documentclass[3p,times,11pt]{elsarticle}
\usepackage[titletoc,title]{appendix}
\usepackage{float}
\usepackage{hyperref}
\usepackage{placeins}
\usepackage{wrapfig,blindtext}
\usepackage{bm}
\usepackage{graphicx} 
\usepackage{epstopdf} 
\usepackage{subfigure} 
\usepackage{enumerate}
\usepackage{pslatex} 
\usepackage{amsmath,amssymb}
\numberwithin{equation}{section}
\usepackage[skip=0pt]{caption}

\usepackage{amsfonts,amsthm} 
\usepackage[english]{babel} 

\usepackage[dvipsnames]{xcolor}

\usepackage{pdfpages}
\usepackage{float}
\usepackage{comment}
\usepackage{nomencl}
\hypersetup{
    colorlinks=true,       
    linkcolor=blue,          
    citecolor=blue,        
    filecolor=magenta,      
    urlcolor=cyan           
}

\usepackage{comment}

\usepackage{nomencl}
  \makenomenclature

\journal{}

\begin{document}

\begin{frontmatter}

\title{\sc A priori estimation of memory effects in coarse-grained nonlinear systems using the Mori-Zwanzig formalism}

\author{Ayoub Gouasmi}
\author{Eric J. Parish}
\author{Karthik Duraisamy}
\address{Department of Aerospace Engineering, University of Michigan, Ann Arbor, USA}




\begin{abstract}
\indent Reduced Order Models (ROMs) of complex, nonlinear dynamical systems often require closure, which is the process of representing the contribution of the unresolved physics on the resolved physics. The Mori-Zwanzig (M-Z) procedure allows one to write down formally closed evolution equations for the resolved physics. In these equations, the unclosed terms are  recast as a memory integral involving the past history of the resolved variables, and a ``noise" term. While the M-Z procedure does not directly reduce the complexity of the original system, these equations can serve as a mathematically consistent starting point to develop closures based on approximations of the memory. In this scenario, a priori knowledge of the memory kernel, which is not explicitly known for nonlinear systems, is of paramount importance to assess the validity of a memory approximation. Unraveling the memory kernel requires the determination of the orthogonal dynamics which is a projected high-dimensional partial differential equation  that is not tractable in general. A method to estimate the memory kernel a priori, using full-order solution snapshots, is proposed. The main idea is to solve a pseudo orthogonal dynamics equation, that has a convenient Liouville form, instead of the original one. This ersatz arises from the assumption that the semi-group of the orthogonal dynamics operator is a composition operator, akin to semi-groups of Liouville operators, for one observable. The method is exact in the linear case where the kernel is known explicitly. Results for under-resolved simulations of the Burgers and Kuramoto-Sivashinsky equations demonstrate that the proposed technique can accurately reconstruct the transfer of information between the resolved and unresolved dynamics through memory, and provide valuable information about the kernel.
\end{abstract}

\end{frontmatter}

\section{Introduction}
Complex dynamical systems are typically prohibitively expensive to solve using direct simulation approaches, insofar as real-world science and engineering applications are concerned. This fact has motivated the development of manifold reduced-order modeling (ROM) techniques based on  Proper Orthogonal Decomposition (POD) \cite{Lumley, Sirovich, Aubry, Berkooz}, Balanced Proper Orthogonal Decomposition (BPOD) \cite{Willcox, Rowley}, balanced truncation \cite{Gugercin, Moore}, the reduced basis method \cite{Veroy, Rozza}, and Krylov subspaces \cite{Beattie}. In projection-based ROM of partial differential equations,  the governing equations are projected onto an alternate space of solutions, in which the state of the system can be represented by an equivalent set of distinct modes.  This typically results in a system of Ordinary Differential Equations (ODE):
\begin{equation}\label{eq:genODE}
\frac{d}{dt}\phi(t) = R(\phi(t)) \ \ ; \ \ \phi(0) = x,
\end{equation}
where $\phi = \phi(x,t) \in \mathbb{R}^{N}$ is the vector of modal coefficients. 
A small subset $\hat{\phi} = (\phi_1, \ \cdots, \phi_m) \in \mathbb{R}^{m}, \ m \ll N$, of this system is then evolved in time. The underlying assumption is that, at every time instant, the essential or dominant dynamics are captured within these $m$ modes. 
The construction of accurate reduced models for nonlinear problems is, however, rendered challenging by the interaction between the resolved and unresolved physics. In many cases, the evolution equation of $\hat{\phi}$ involves the discarded modes $\tilde{\phi} = (\phi_{m+1}, \ \cdots, \phi_N) \in \mathbb{R}^{N-m}$. This contribution has to be modeled, and the process is referred to as closure. This is an outstanding challenge that reduced-fidelity approaches such as Large Eddy Simulation (LES) \cite{Smagorinsky, Deardorff} must also deal with.

The present work is inspired from the optimal prediction framework developed by the Chorin group \cite{Chorin, Chorin2, Chorin3, Chorin4}. This framework is a reformulation of the Mori-Zwanzig (M-Z) formalism  \cite{Mori, Zwanzig} of non-equilibrium statistical mechanics. The M-Z formalism provides a mathematically consistent technique for the development of reduced-order models for dynamical systems. A key element is the use of projection operators $\mathcal{P}$ in phase space. Functions $f(x)$ of the initial state are projected onto a subspace of functions of the resolved part $\hat{\phi}(0) = \hat{x}$ only. The formalism consists of separating the state variables into a resolved and unresolved set. The nonlinear dynamical system is recast as a Liouville equation; a linear partial differential equation (PDE) in the space of the initial conditions. An exact system, known as the generalized Langevin equation, can then be derived that describes the evolution of the resolved variables as a function of only the resolved variables. In this setting, the effect of the unresolved modes on the resolved modes appears as a non-local memory integral and a noise term. The M-Z approach - by itself - does not lead to a reduction in computational complexity as the evaluation of the memory and noise terms requires the solution of the orthogonal dynamics equation.
The formalism, instead, provides a starting point for the development of closure models.

It is important to emphasize that, while this work is inspired from the Chorin framework, and the Mori-Zwanzig formalism by transitivity, the authors are considering applications to a different, simpler in concept, type of problem. The Mori-Zwanzig procedure was developed in the non-equilibrium statistical mechanics community, where the goal is to solve for probability density functions and time correlation functions of non-equilibrium systems \cite{Zwanzig0}. Originally, the procedure was limited to Hamiltonian dynamical systems.
It was only when Chorin extended this result to general time-dependent systems (\cite{Chorin}, chapter 9) that extensions to non-hamiltonian problems such as the viscous Burgers equation \cite{Bernstein}, the Kuramoto-Sivashinsky \cite{KS} equation or the Euler equations \cite{Stinis2} were considered. Chorin's framework was developed for optimal prediction, that is estimating the solution of nonlinear time-dependent problems such as eq. (\ref{eq:genODE}) where a full-order solution is not affordable and the unresolved part $\tilde{x}$ of the initial conditions is uncertain. Chorin and co-workers  developed a method to compute mean solutions with respect to initial probability distributions. The present work considers nonlinear, non-Hamiltonian, time-dependent systems where the initial conditions are fully resolved, i.e. $\tilde{\phi}(0) = \tilde{x} = 0$, with no uncertainty.


In optimal prediction, the projection operator is a conditional expectation $(\mathcal{P}f)(x) = E [ f(x) | \hat{x}]$. It can also be approximated by a finite rank projection. For the problem that the present work is concerned with,  the projector is a simple truncation $(\mathcal{P}f)(\hat{x},\tilde{x}) = f(\hat{x},0)$. It has been used in \cite{Parish, StinisP}. Although the M-Z approach is independent of the projector, its viability is dependent upon the well-posedness of the orthogonal dynamics equation \cite{Givon}, which depends on the projector. Givon, Hald and Kupfermann \cite{Givon} studied the existence of solutions to the orthogonal dynamics for Hamiltonian systems. They proved the existence of classical solution for finite-rank projections and the existence of weak solutions for the conditional expectation. Even though the truncation projector can be seen as a limiting case of the conditional expectation, non-Hamiltonian systems are considered in this work. 

Solution techniques for the orthogonal dynamics have been developed for finite-rank projection operators. One of the most notable has been devised by Chorin \cite{Chorin}. It involves the decomposition of the solution using a finite-rank orthonormal basis and the evolution in time of a set of Volterra equations for the basis coefficients. One of the main issues with this approach lies in the high expense associated with the inner product and the number of basis functions required for the decomposition to be representative. Twenty-one basis functions were used for the 4-unknown Hamiltonian system studied in \cite{Chorin}. For the Burgers equation, Bernstein \cite{Bernstein} showed that unless the number of resolved modes is very small, a similar approach for the truncation projection operator case would not be computationally tractable. The number of basis functions needed would typically escalate for higher-dimensional problems. In the Molecular Dynamics community, Darve \cite{Darve} developed a discrete version of the Mori-Zwanzig formalism and applied it to compute the kernel in GLEs of Hamiltonian systems. This methodology was developed for the conditional expectation projector, and it was shown that the solution to the discrete orthogonal dynamics is the same as that of the continuous version in the asymptotic limit. He also points out the high computational expense if the dimension of the resolved space is large. 

At this juncture, it is important to emphasize that the Mori-Zwanzig formalism does not lead to a ROM directly, as it does not address the issue of how to evaluate the outstanding memory term in a tractable fashion, nor whether it is even possible. A priori knowledge of the memory kernel can help answer these two questions. This is the primary purpose of this work. A method is developed to numerically estimate the memory kernel a priori for nonlinear systems. The key idea is to consider an ersatz of the orthogonal dynamics equations that can be solved using the method of characteristics. This ersatz arises from the assumption that the semi-group generated by the orthogonal dynamics is a composition operator. Furthermore, it will be shown that evaluating the memory kernel amounts to computing sensitivities of the solution to the orthogonal dynamics with respect to the initial conditions in a direction that is imposed by the ROM solution. 

The remainder of this paper will be organized as follows: Section 2 will introduce the the M-Z formalism. In Section 3, an approach to estimate the memory kernel  will be presented. In Section 4, our procedure will be applied to the Brusselator, the Viscous Burgers equation and Kuramoto-Sivashinsky (K-S). Conclusions and perspectives will be provided in Section 5.

\section{The Mori-Zwanzig Formalism}
\subsection{Model reduction and memory}
\begin{figure}[h!]
    \centering
    \subfigure[t = 2 s]{\includegraphics[width = .48\linewidth, height = 0.24\textheight]{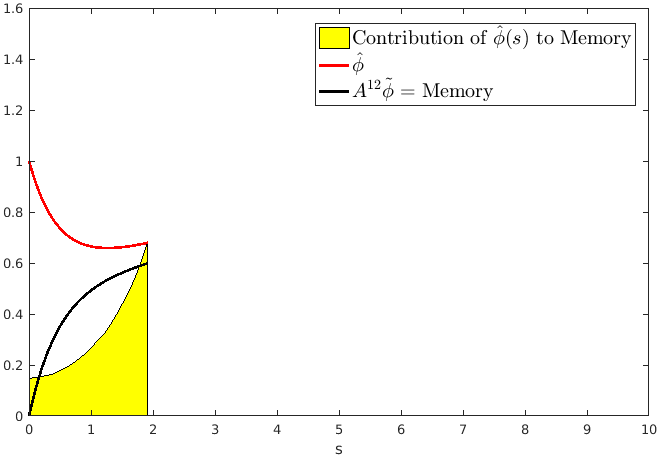}}
    \subfigure[t = 6 s]{\includegraphics[width = .48\linewidth, height = 0.24\textheight]{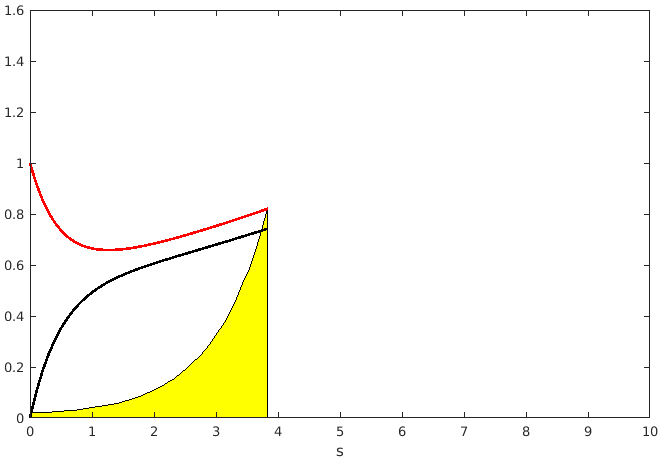}}
    \subfigure[t = 10 s]{\includegraphics[width = .48\linewidth, height = 0.24\textheight]{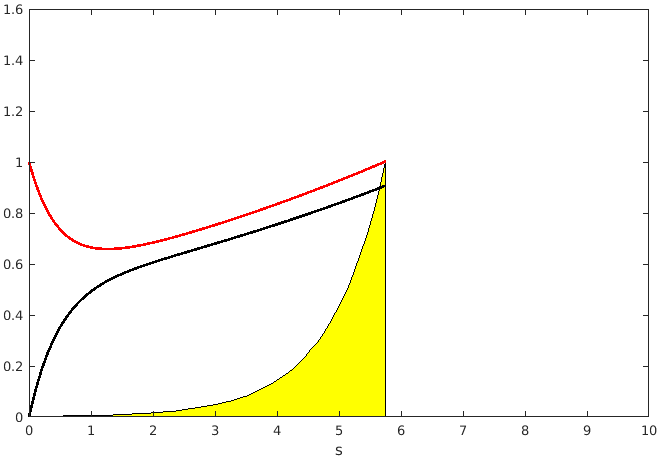}}
    \subfigure[t = 16 s]{\includegraphics[width = .48\linewidth, height = 0.24\textheight]{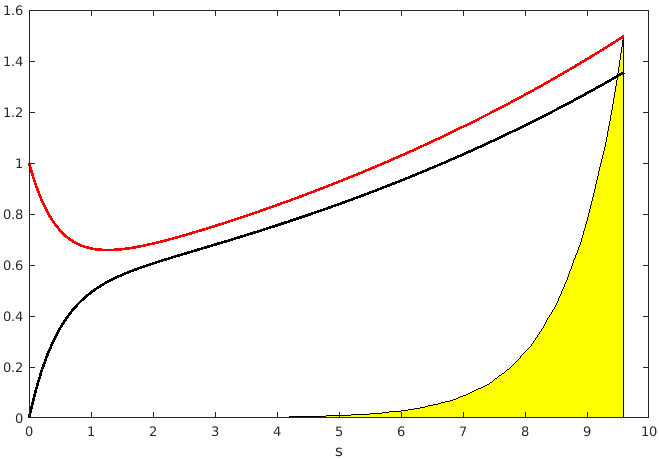}}
    \caption{Linear system: Snapshots of the memory kernel. $\forall s \in [0, \ t]$, the contribution of $\hat{\phi}(s)$ to the memory at $t$ is exactly $A^{12}e^{A^{22}(t-s)}A^{21}\hat{\phi}(s)$. Memory appears finite, and the decay profile of the kernel suggests that the memory can be approximated using the contribution of $\hat{\phi}(t)$ and a measure of $\tau$ only.}
    \label{fig:Line_decay_hist}
\end{figure}

As an introductory example \cite{Zwanzig0}, consider  a linear dynamical system:
\begin{equation}\label{eq:lin}
\begin{cases}
\begin{gathered}
\frac{d}{dt}\hat{\phi}   = A^{11}\ \hat{\phi} \ + \ A^{12} \ \tilde{\phi}, \\ 
\frac{d}{dt}\tilde{\phi} = A^{21}\ \hat{\phi} \ + \ A^{22} \ \tilde{\phi}, \\
\phi(0) = x,
\end{gathered}
\end{cases}
\end{equation}
with $A^{11} \in \mathbb{R}^{m \times m}$, $A^{12} \in \mathbb{R}^{m \times (N -m)}$, $A^{21} \in \mathbb{R}^{(N - m) \times m}$ and $A^{22} \in \mathbb{R}^{(N-m) \times (N-m)}$. To solve directly for $\hat{\phi}$, one has to model the "unclosed" term $A^{12} \ \tilde{\phi}$ as a function of $\hat{\phi}$ only. This can be accomplished by first solving the equation for $\tilde{\phi}$ assuming knowledge of $\hat{\phi}$:
\begin{equation}\label{eq:DysonLin}
\tilde{\phi}(t) = e^{\ A^{22}t}\tilde{x} \ + \int_0^t e^{\ A^{22}(t - s)} A^{21} \hat{\phi}(s) ds.
\end{equation}
By injecting eq. (\ref{eq:DysonLin}) in the closure term, one obtains: 
\begin{equation} \label{eq:MZLin}
\frac{d}{dt}\hat{\phi} = \ A^{11} \ \hat{\phi} \ + \  \int_0^t A^{12} e^{\ (t-s)A^{22}} A^{21} \hat{\phi}(s) ds \ + \ A^{12} e^{\ A^{22}t}\tilde{x}.
\end{equation}
Three distinct terms feature in this equation: a Markovian term, a memory involving the past history of the resolved physics and a term that encapsulates the influence of the initial conditions. In an optimal prediction context, the latter term is interpreted as noise produced by the uncertainty in the initial conditions. In our context, there is no noise contribution because the initial condition is fully resolved $\tilde{x} = 0$. The closure term is exclusively memory. Eq. (\ref{eq:MZLin}) is often called a Generalized Langevin equation (GLE). \\
\indent In eq. (\ref{eq:MZLin}), the dependency on $\tilde{\phi}$ has been removed without introducing any approximation. However, eq. (\ref{eq:MZLin}) is not a ROM because of the inherent cost. It requires multiple evaluations of the exponential of a (large) $(N - m)\times (N - m)$ matrix multiplied by a $(N - m)\times m$ matrix. The cost of evolving this ROM in time should be comparable to that of computing $A^{11} \ \hat{\phi}$. \\ 
\indent To make eq. (\ref{eq:MZLin}) a genuine ROM, the memory term has to be approximated. Indeed, the approximation depends on the integrand $A^{12}e^{\ (t-s)A^{22}} A^{21} \hat{\phi}(s)$, which is called the memory kernel. Let's assume that  the entire set of eigenvalues of $A^{22}$ is strictly negative. Then a relevant assumption to make is that the memory has finite support, i.e. there exists some memory length $\tau > 0$ such that:
\begin{equation*}
    \int_0^t A^{12} e^{\ (t-s)A^{22}} A^{21} \hat{\phi}(s) ds \ \approx \ \int_{t-\tau}^t A^{12} e^{\ (t-s)A^{22}} A^{21} \hat{\phi}(s) ds.
\end{equation*}
\indent This is highlighted in figure \ref{fig:Line_decay_hist} for matrices $(A^{11},\ A^{12},\ A^{21},\ A^{22}) = (-0.8, 1, 1, -1)$. The initial conditions are $(\hat{x},\tilde{x}) = (1,0)$. The subgrid term $A^{12} \tilde{\phi}$ at time t is equal to the area of the shaded yellow region $A^{12}e^{(t-s)A^{22}}A^{21}\hat{\phi}(s)$ over the elapsed time range $s \in [0, \ t]$. In these snapshots, the memory does appear to be finite. In addition, the decay profile of the kernel, from its most recent ($s = t$) to its oldest ($s = 0$) contributions, suggest that the subgrid/memory term could be reasonably approximated using the most recent kernel contribution and a good guess of the memory length $\tau$. In other words, if the kernel has a simple decay profile (exponential, for instance) one could consider a closure model of the form:
\begin{equation*}
    A^{12} \tilde{\phi}(t) = \int_0^t A^{12} e^{\ (t-s)A^{22}} A^{21} \hat{\phi}(s) ds \ \approx \ A^{12} \ C(\tau) \ A^{21} \hat{\phi}(t),
\end{equation*}
with $C(\tau)$ being either a scalar or a $(N - m) \times (N - m)$ matrix, depending on the level of sophistication desired. The memory length $\tau$ does not have to be same for all $m$ components of the memory kernel, nor does it have to be constant throughout the run. If $C(\tau)$ is simple enough the evaluation of the memory approximation can be made reasonable enough for use as a ROM closure model. \\  
\indent It is important to appreciate the fact that the above discussion on developing closure models based on memory was straightforward only because rewriting eq. (\ref{eq:lin}) as a GLE, with all outstanding terms explicitly known, was possible. As one can anticipate, however, for nonlinear systems, additional challenges arise. \\ 
\indent A first nonlinear ODE considered in this work is the Brusselator ~\cite{Prigogine}: 
\begin{equation}\label{eq:Bru}
\begin{cases}
\begin{aligned}
\frac{d}{dt}\phi_1 =& \ A \ - \ B \ \phi_1 \ - \ \phi_1 \ + \  \phi_1^2 \ \phi_2, \\
\frac{d}{dt} \phi_2 =& \ B \ \phi_1 \ - \ \phi_1^2 \ \phi_2,
\end{aligned}
\end{cases}
\end{equation}
where $\hat{\phi} = \phi_1$, $\tilde{\phi} = \phi_2$, $m = 1$ and $N= 2$. $A$ and $B$ are constants. The Brusselator will serve as a toy problem to familiarize the reader with the concepts and the method presented in this work. \\
\indent It is possible to write a GLE for the evolution of $\hat{\phi}$ for nonlinear systems. For instance, it is possible to rewrite the subgrid term $\hat{\phi}^2 \ \tilde{\phi}$ in the Brusselator as:
\begin{equation*}
    \hat{\phi}(t)^2 \ \tilde{\phi}(t) = \int_0^{t} K_1(\hat{\phi}(s), t - s) ds + F_1(x,t),
\end{equation*}
with $F_1(x,t)$ dropping if $\tilde{x} = 0$. Further details are provided in the next sections. The technique used to obtain the above GLE is based on projection operator methods \cite{Hynes} and was derived by Mori \cite{Mori} and \cite{Zwanzig0} for non-equilibrium, Hamiltonian systems. Chorin \cite{Chorin} extended it to general dynamical systems. In this work, important concepts of the M-Z formalism are discussed before presenting the derivation \cite{Chorin, Chorin3} of the GLE in the general case. \\
\indent Contrary to the linear case, very little is known about the memory kernel (contribution of $\hat{\phi}(t)$ only) and the noise (value at $t = 0$ only). More details are provided in the following sections. In short, memory approximations are hard to devise in this case because of the outstanding challenge posed by an auxiliary orthogonal dynamics equation. The goal of the present work is to develop a numerical method to estimate the memory kernel for nonlinear systems.  

\subsection{Liouville equivalence}
\indent The Liouville operator associated with eq. (\ref{eq:genODE}) is
\begin{equation*}
    \mathcal{L} = \sum_{i=1}^N R_i(x) \frac{\partial}{\partial x_i}.
\end{equation*}
Consider the following linear PDE in phase space:
\begin{align}\label{eq:genLiouville}
\begin{split}
\frac{\partial}{\partial t} u(x,t) = \mathcal{L} u(x,t) \ \ ; \ \  
u(x,0) = g(x),
\end{split}
\end{align}
where $g(x)$ is a function of the initial conditions. Eq. (\ref{eq:genLiouville}) is the Liouville equation. Chorin showed that, for any regular function g(x), the solutions of equations (\ref{eq:genODE}) and (\ref{eq:genLiouville}) satisfy: 
\begin{equation}\label{eq:LEQ}
    u(x,t) = g(\phi(x,t)).
\end{equation}
We refer to this result as the Liouville equivalence. In the case $g(x) = x$, the solution to the Liouville equation is the same as that of eq. (\ref{eq:genODE}). The Liouville equivalence demonstrates that the solutions $u(x,t)$ of eq. (\ref{eq:genLiouville}) for all functions $g(x)$ are known once the solution $\phi(x,t)$ to eq. (\ref{eq:genODE}) is known. This is an instance of the method of characteristics for advection equations. This equivalence is a fundamental step in deriving the GLE (section 2.4).\\
\indent In semi-group notation, the solution of the Liouville equation can be rewritten as:
\begin{equation*}
u(x,t) = e^{t\mathcal{L}} u(x,0) = e^{t\mathcal{L}} g(x),
\end{equation*}
where $e^{t\mathcal{L}}$ is the semi-group generated by $\mathcal{L}$. The Liouville equivalence endows this semi-group with a remarkable composition property: 
\begin{equation}\label{eq:Lfeat}
e^{t\mathcal{L}} g(x) = g(e^{t\mathcal{L}}x).
\end{equation}
This property is of central interest to this work. \\
\indent As a side note, the semi-group $e^{t\mathcal{L}}$ is also known as the Koopman semi-group \cite{Koopman1, Koopman2} associated with the continuous-time dynamical system (\ref{eq:genODE}). 
\subsection{Projection Operators}
Projection operators \cite{Hynes} are central in non-equilibrium problems and the Mori-Zwanzig formalism. For purposes of introduction, the optimal prediction setting of Chorin \cite{Chorin} will be initially assumed. \\
\indent Consider the initial conditions x to be random variables drawn from some probability distribution $\mu(x)$ and a probability density function $\rho(x)$. Denote $\Gamma$, the vector space in which $x$ and $\phi(t)$ lie ($\mathbb{R}^N$ in our case). For a function g(x) that operates on $\Gamma$, define the expected value of g by:
\begin{equation*}
    E[g] = \int_{\Gamma} g(x) \rho(x) dx.
\end{equation*}
Denote $L^2$ the Hilbert space of functions on $\Gamma$ endowed with the inner product $(f, g) = E[fg]$. In the Mori-Zwanzig formalism, a projector $\mathcal{P}$ that transforms functions in $L^2$ into functions of the resolved component $\hat{x} \in \mathbb{R}^{m}$, is considered. There are several different projectors \cite{Chorin} for which the M-Z formalism has been considered:
\begin{enumerate}[(I)]
    \item In irreversible statistical mechanics, the linear projection is given by:
    \begin{equation*}
    \forall f \in L^2, (\mathcal{P}f)(\hat{x}) = \sum_{i, j = 1}^m (x_i, x_j)^{-1} (f, x_i) x_i.
    \end{equation*}
    \item In optimal prediction, the conditional expectation of $f$ given $\hat{x}$, is given by:
    \begin{equation*}
        \forall f \in L^2, \ (\mathcal{P}f)(\hat{x}) = E[f|\hat{x}] = \frac{\int f(x)\rho(x)d\tilde{x}}{\int \rho(x)d\tilde{x}}, \ d\tilde{x} = dx_{m+1}\cdots dx_N.
    \end{equation*}
    This is also the orthogonal projection onto the span of all functions of $\hat{x}$.
    \item The conditional expectation can be approximated by a finite-rank projection onto the span of an orthonormal set of function $h_{k}(\hat{x}), k = 1 \dots K$:
    \begin{equation*}
        \forall f \in L^2, \ (\mathcal{P}f)(\hat{x}) = \sum_{k = 1}^K (f, h_{\nu}) h_{k}(\hat{x}).
    \end{equation*}
\end{enumerate}
Define $\mathcal{Q} = \mathcal{I} - \mathcal{P}$ the complimentary projector, where $\mathcal{I}$ is the identity operator. \\
\indent In this work, we focus on under-resolved simulations of nonlinear time-dependent systems with fully-resolved initial conditions (i.e. $\tilde{x} = 0$). The projector to be used for that setup is a truncation:
\begin{equation} \label{eq:P}
(\mathcal{P}f)(\hat{x}) =  f(\hat{x},0) =  f(\boldsymbol{\hat{x}}).
\end{equation}
This projector can be seen as a conditional expectation in the limit where the probability density function is a delta function centered at $\tilde{x} = 0$.

\subsection{Generalized Langevin Equation}
The Liouville eq. (\ref{eq:genLiouville}) for $g(x) = x_j$ is given by:
\begin{equation*}
    \frac{\partial}{\partial t} e^{t\mathcal{L}}x_j = \mathcal{L} e^{t\mathcal{L}}x_j = e^{t\mathcal{L}}\mathcal{L}x_j.
\end{equation*}
The right-hand side term of the above equation can decomposed using $\mathcal{I} = \mathcal{P} + \mathcal{Q}$:
\begin{equation}\label{eq:L}
    \frac{\partial}{\partial t} e^{t\mathcal{L}}x_j = e^{t\mathcal{L}}\mathcal{PL}x_j + 
       e^{t\mathcal{L}}\mathcal{QL}x_j.
\end{equation}
The first term $e^{t\mathcal{L}}\mathcal{PL}x_j$ is a function of $\hat{\phi}$ only. The point of recasting the dynamical system (\ref{eq:genODE}) as a linear PDE is that, using the Dyson formula \cite{Evans}:
\begin{equation*}
    e^{t\mathcal{L}} = e^{t\mathcal{QL}} + \int_0^t e^{(t-s)\mathcal{L}}\mathcal{PL}e^{s\mathcal{QL}}\ ds,
\end{equation*}
the second term $e^{t\mathcal{L}}\mathcal{QL}x_j$ in eq. (\ref{eq:L}) that depends on the full solution $\phi \in \mathbb{R}^N$ can be rewritten as:
\begin{equation*}
    e^{t\mathcal{L}}\mathcal{QL}x_j = e^{t\mathcal{QL}}\mathcal{QL}x_j + \int_0^t e^{(t-s)\mathcal{L}}\mathcal{PL}e^{s\mathcal{QL}}\mathcal{QL}x_j \ ds.
\end{equation*}
The first term $F_j(x,t) = e^{t\mathcal{QL}}\mathcal{QL}x_j$ models the influence of the initial conditions. It is often interpreted as the noise produced by the uncertainty in the initial conditions. The second term is the general expression for the memory, which involves the past history of resolved physics $\hat{\phi}(s)$ only. Ultimately, the dynamical system (\ref{eq:genODE}) has been recast as a set of GLEs for the state components $\phi_j(x, t) = e^{t \mathcal{L}}x_j$:  
\begin{equation}\label{eq:GLE}
\frac{\partial}{\partial t} e^{t\mathcal{L}}x_j = e^{t\mathcal{L}}\mathcal{PL}x_j \ + \ \int_0^t e^{(t - s)\mathcal{L}}\mathcal{PL}e^{s \mathcal{QL}}\mathcal{QL}x_j \ ds \ + \ e^{t\mathcal{QL}}\mathcal{QL}x_j,
\end{equation}
or, in a more compact form:
\begin{gather*}
\frac{\partial}{\partial t} \phi_j(x,t) = \mbox{R}_j(\hat{\phi}(x,t)) \ + \ \int_0^t K_j(\hat{\phi}(x,s), t - s) \ ds \ + \ F_j(x, t), \\
\mbox{R}_j(\hat{x}) = (\mathcal{P}R_j) (\hat{x}), \ F_j(x,t) = e^{t\mathcal{QL}}\mathcal{QL}x_j, \ K_j(\hat{x},t) = \mathcal{PL}F_j(x,t).
\end{gather*}
In the RHS, the first term is the Markovian contribution that depends on $\hat{\phi}(x,t)$ only. The second term is the memory that will involve contributions from all past values of the resolved physics $\hat{\phi}(x,s), s \in [0, \ t]$ only. The relationship between the memory kernel $K_j$ and the noise $F_j$ is referred to as a "fluctuation dissipation theorem" \cite{Chorin, Chorin2, Chorin3, Zwanzig0, Stuart} in the original context of statistical mechanics. \\
\indent The interpretation and implementation of each term is imposed by the projector $\mathcal{P}$ \cite{Stuart}. The projector is chosen such that $\mathcal{P}\hat{\phi}(x, t)$ is the target quantity the user wants to solve for. The evolution equation for the target is obtained by projecting the GLE:
\begin{equation}
\frac{\partial}{\partial t} \mathcal{P} e^{t\mathcal{L}}x_j = \mathcal{P}e^{t\mathcal{L}}\mathcal{PL}x_j \ + \ \mathcal{P} \int_0^t e^{(t - s)\mathcal{L}}\mathcal{PL}e^{s \mathcal{QL}}\mathcal{QL}x_j \ ds.
\end{equation}
As is often the case in large eddy simulations of fluid flow, the authors choose  $\mathcal{P}\hat{\phi}(x,t) = \hat{\phi}(\hat{x},t)$. The corresponding equation is therefore: 
\begin{gather}\label{eq:PM-ZEvo}
\frac{\partial}{\partial t} \phi_j(\hat{x},t) = R_j(\hat{\phi}(\hat{x},t)) + \int_0^t K_j(\hat{\phi}(\hat{x},s),t-s) ds.
\end{gather}

\indent 
As explained in section (1.1) for the linear case, eq. (\ref{eq:PM-ZEvo}) is not a reduced order model. This is because of the intrinsically high cost associated with the evaluation of the kernel (despite being a function of $\hat{\phi} \in \mathbb{R}^{m}$ only!). For instance, evaluation of the kernel in the linear case involved operating on matrices of size much larger than $m$. In the nonlinear case, evaluating the kernel requires the solution of the orthogonal dynamics equation.                                                       

\subsection{Orthogonal Dynamics}
$F(x,t) = (\ F_j(x,t)\ )_{1 \leq j \leq N} \in \mathbb{R}^{N}$ is the solution of  a linear PDE:
\begin{equation} \label{eq:OrtDyna}
\frac{\partial}{\partial t}F(x,t) = \mathcal{QL}F(x,t) \ \ ; \ \
F(x,0) = \mathcal{QL}x.
\end{equation}
 This equation is called the orthogonal dynamics because its solution lies in the orthogonal space of $\mathcal{P}$ at all instants. 

\indent Givon, Hald and Kupfermann \cite{Givon} studied the existence of solutions to the orthogonal dynamics for Hamiltonian systems. They proved the existence of classical solution for finite-rank projections and the existence of weak solutions for the conditional expectation. Even though the truncation projector can be seen as a limiting case of the conditional expectation, non-Hamiltonian systems are considered in this work. Therefore existence of solutions to the associated orthogonal dynamics cannot be supported using their work. \\
\indent It is possible to find a solution to the orthogonal dynamics when the dynamical system is linear. Consider eq. (\ref{eq:lin}), and define matrices $A$ and $\tilde{A}$ $\in \mathbb{R}^{N \times N}$ as:
\begin{equation*}
A = 
\begin{bmatrix}
A^{11} & A^{12} \\ A^{21} & A^{22}
\end{bmatrix}, \
\tilde{A} = 
\begin{bmatrix}
0 & A^{12} \\ 0 & A^{22}
\end{bmatrix}. 
\end{equation*}
The Liouville operator is given by:
\begin{equation*}
    \mathcal{L} = \sum_{i = 1}^N \big (\sum_{j = 1}^N A_{ij}x_j \big) \frac{\partial}{\partial x_i}
\end{equation*}
\indent We note that:
\begin{equation*}
    \mathcal{QL}x = R(x) - \mathcal{P}R(x) = A \begin{bmatrix} 0 \\ \tilde{x} \end{bmatrix} = \tilde{A} \ x.
\end{equation*}
therefore:
\begin{equation*}
    \forall n, \ (\mathcal{QL})^n x = \tilde{A}^n x,
\end{equation*}
and it can be shown, by induction, that:
\begin{equation*}
    \tilde{A}^n = \begin{bmatrix}
0 & A^{12} (A^{22})^{n-1} \\ 0 & (A^{22})^n
\end{bmatrix}.
\end{equation*}
For any matrix $B$, the series $\sum_{n = 0}^{\infty} \frac{t^n}{n!}B^n$ converges to $e^{t B}$. Therefore, so does the series $\sum_{n = 1}^{\infty}\frac{t^n}{n!}(\mathcal{QL})^n x$ and:
\begin{equation*}
    F(x,t) = \sum_{n = 1}^{\infty}\frac{t^n}{n!}(\mathcal{QL})^n x= \begin{bmatrix}
0 & A^{12} e^{t A^{22}} \\ 0 & A^{22} e^{t A^{22}}
\end{bmatrix}x = \begin{bmatrix} A^{12} \\ A^{22}\end{bmatrix} e^{t A^{22}} \tilde{x},
\end{equation*}
is a solution. Then:
\begin{align*}
\mathcal{L}F(x,s) =& \   \begin{bmatrix} A^{12} \\ A^{22}\end{bmatrix} e^{ \ A^{22}s}  \sum_{i = 1}^N \big (\sum_{j = 1}^N A_{ij} x_j\big) \frac{\partial}{\partial x_i}  \tilde{x}
= \begin{bmatrix} A^{12} \\ A^{22}\end{bmatrix} e^{ \ A^{22}s}  \sum_{i = m+1}^N \big (\sum_{j = 1}^N A_{ij}x_j \big) \frac{\partial}{\partial x_i}  \tilde{x} = \begin{bmatrix} A^{12} \\ A^{22}\end{bmatrix} e^{ \ A^{22}s} [A^{21} \ A^{22}]x
\end{align*}
Applying the projector $\mathcal{P}$, one gets:
\begin{equation*}
K(\hat{x},s) = \mathcal{PL}F(x,s) =  \begin{bmatrix} A^{12} \\ A^{22}\end{bmatrix} e^{ \ A^{22}s} [A^{21} \ A^{22}] 
\begin{bmatrix}
 \hat{x} \\ 0_{(N-m) \times 1}\end{bmatrix} = \begin{bmatrix} A^{12} \\ A^{22}\end{bmatrix} e^{ \ A^{22}s} A^{21} \hat{x}.
\end{equation*}
Ultimately,
\begin{equation*}
K(\hat{\phi}(s),t-s) = \begin{bmatrix} A^{12} \\ A^{22}\end{bmatrix} e^{ \ A^{22}(t - s)} A^{21} \hat{\phi}(s).
\end{equation*}
is an exact expression. The first $m$ components correspond to what is in eq. (\ref{eq:MZLin}). \\  
\indent As a side note, the above solution method did not assume analyticity of $e^{t\mathcal{QL}}$, meaning that this semi-group could be described as:
\begin{equation*}
    e^{t\mathcal{QL}} = \sum_{n = 0}^{\infty}\frac{t^n}{n!}(\mathcal{QL})^n.
\end{equation*}
In semi-group theory \cite{Pazy}, such an expansion is convergent if the generator $\mathcal{QL}$ is a bounded operator on a Banach space.  To find a solution for the linear case, we just showed that the series $\frac{t^n}{n!}(\mathcal{QL})^n x$  was convergent. 
\section{A Priori Estimation of the Memory Kernel}
This work seeks to further an approach Parish and  Duraisamy \cite{Parish} hinted at to solve the orthogonal dynamics and suggests a procedure to estimate the memory kernel in nonlinear dynamical systems. \\
\indent An equivalent form of eq. (\ref{eq:PM-ZEvo}) that is more fitting to the authors' viewpoint is the following:
\begin{gather}\label{eq:MZref}
w_j = R_j(\phi(\hat{x}, t)) - R_j(\hat{\phi}(\hat{x},t)) = \int_0^t K_j(\hat{\phi}(\hat{x},s),t-s) ds.
\end{gather}
The LHS term is the general expression for the subgrid terms, denoted $w_j$, that require modeling in ROMs for $j = 1 \dots m$. The Mori-Zwanzig formalism shows that it can be exactly represented as a memory integral. As a reminder, the noise term dropped because we are within the setting the projector $\mathcal{P}$ imposes, which is fully resolved initial conditions. \\
\indent The goal of this work is not to derive closure models based on memory approximations as in \cite{Bernstein, Parish, Chorin, KS}. No closure model has been derived or tested in ROM simulations. Full-order simulations have been run to provide for the "exact" ROM solution snapshots $\hat{\phi}(\hat{x},t)$ and the "exact" subgrid terms $w_j$. The  focus of this work is in estimating the memory kernel. The goal is to extract the kernel as in section (2.1), but for nonlinear systems. \\
\indent The memory kernel is unknown for nonlinear systems. Eq. (\ref{eq:MZref}) is therefore the only verification model available. In the numerical results section, the exact ROM snapshots and subgrid terms are first computed. Then the procedure is applied to estimate the memory kernel $K_j(\hat{\phi}(\hat{x},s), t - s)$ at the required points and time ranges. 

\subsection{Pseudo orthogonal dynamics}
The solution of the orthogonal dynamics can be formally expressed as:
\begin{equation*}
F(x,t) = e^{t\mathcal{QL}}F(x,0),
\end{equation*}
where
\begin{equation*}
F(x,0) = \mathcal{QL}x = \mathcal{Q}R(x) = R(x) - \mathcal{P}R(x) = R(x) - R(\hat{x}).
\end{equation*}
Earlier, the composition property of the semi group $e^{t\mathcal{L}}$ was introduced:
\begin{equation*}
e^{t\mathcal{L}}g(x) = g(e^{t\mathcal{L}}x) = g(\phi(x,t)).
\end{equation*}
The main assumption of this work consists in granting the semi-group of the orthogonal dynamics  $e^{t\mathcal{QL}}$ the composition property for $g(x) = \mathcal{QL}x$. Defining $\phi^{\mathcal{Q}}(x,t) =  e^{t\mathcal{QL}}x \in \mathbb{R}^{N}$, one has:
\begin{equation*}
F(x,t) = e^{t\mathcal{QL}}F(x,0) = F(e^{t\mathcal{QL}}x,0) = F(\phi^{\mathcal{Q}}(x,t),0) = R(\phi^{\mathcal{Q}}(x,t)) - R(\hat{\phi}^{\mathcal{Q}}(x,t)),
\end{equation*}
where $\hat{\phi}^{\mathcal{Q}}$ is defined in the same manner as $\hat{\phi}$. By definition, $\phi^{\mathcal{Q}}(x,t)$ satisfies:
\begin{equation}\label{eq:OrtElem}
\frac{\partial}{\partial t} \phi^{\mathcal{Q}} = \mathcal{QL} \phi^{\mathcal{Q}} \ \ ; \ \ \phi^{\mathcal{Q}}(0) = x.
\end{equation}
Using the expression $\phi^{\mathcal{Q}}(x,t) = e^{t\mathcal{QL}} x$ and the fact that $e^{t\mathcal{QL}}$ and $\mathcal{QL}$ commute, one obtains:
\begin{equation*}
\frac{\partial}{\partial t} \phi^{\mathcal{Q}} = e^{t\mathcal{QL}}\mathcal{QL}x = e^{t\mathcal{QL}}\mathcal{Q}R(x) =  e^{t\mathcal{QL}} (R(x) - R(\hat{x})).
\end{equation*} 
Using composition again, we finally obtain:
\begin{equation}\label{eq:ortODE}
\frac{d}{dt} \phi^{\mathcal{Q}}(t) = R(\phi^{\mathcal{Q}}(t)) - R(\hat{\phi}^{\mathcal{Q}}(t)) \ \ ; \ \ \phi^{\mathcal{Q}}(0) = x.
\end{equation} 
This is an ODE, where the right-hand side (RHS) term is $F(x,t)$. Accordingly, $F_j(x,t)$ can be computed by evolving in time eq. (\ref{eq:ortODE}) and taking the j-th component of the RHS.   
\subsubsection*{Interpretation}
Eq. (\ref{eq:ortODE}) results from granting the semi-group of the orthogonal dynamics a composition property. Let's introduce the Liouville operator $\mathcal{L}^{\mathcal{Q}}$ associated with that ODE:
\begin{equation*}
    \mathcal{L}^{\mathcal{Q}} = \sum_{i = 1}^N \big(R_i(x) - R_i(\hat{x})\big) \frac{\partial}{\partial x_i} = \sum_{i = 1}^N (\mathcal{Q}R_i)(x)  \frac{\partial}{\partial x_i}.
\end{equation*}
Using the Liouville equivalence for $g(x) = F(x,0) = \mathcal{QL}x$, we show that $g(\phi^{\mathcal{Q}}(x,t)) = F(x,t) $ satisfies: 
\begin{align}\label{eq:Ersatz}
\begin{split}
\frac{\partial}{\partial t} F(x,t) = \mathcal{L}^{\mathcal{Q}} F(x,t) \ \ ; \ \  
F(x,0) = \mathcal{QL}x.
\end{split}
\end{align}
Basically, the orthogonal dynamics equation has been replaced with an ersatz (or "pseudo orthogonal dynamics") whose solution can be related to that of eq. (\ref{eq:ortODE}) - that we call pseudo orthogonal ODE - via the Liouville equivalence. \\
\indent Whether or not granting $e^{t\mathcal{QL}}$ a composition property is valid depends on how close the solutions to the orthogonal dynamics eq. (\ref{eq:OrtDyna}) and its ersatz eq. (\ref{eq:Ersatz}) are. For linear systems, the pseudo orthogonal ODE is the following:
\begin{equation*}
\begin{cases}
\begin{gathered}
\frac{d}{dt}\hat{\phi}^{\mathcal{Q}}   = A^{12} \ \tilde{\phi}^{\mathcal{Q}}, \\ 
\frac{d}{dt}\tilde{\phi}^{\mathcal{Q}} = A^{22} \ \tilde{\phi}^{\mathcal{Q}}, \\
\phi^{\mathcal{Q}}(0) = x.
\end{gathered}
\end{cases}
\end{equation*}
and the ersatz solution is given by the RHS, that is:
\begin{equation*}
    F(x,s) = \begin{bmatrix} A^{12} \\ A^{22}\end{bmatrix} \tilde{\phi}^{\mathcal{Q}} = \begin{bmatrix} A^{12} \\ A^{22}\end{bmatrix} e^{t A^{22}} \tilde{x}.
\end{equation*}
This is the same solution that we found in section (2.5). It has already been shown that this solution leads to the exact kernel. \\
\indent For nonlinear systems, however,  solutions of the orthogonal dynamics are unknown. In the absence of reference nonlinear problems where the solution of the orthogonal dynamics is known, an indirect verification approach driven by eq. (\ref{eq:MZref}) is followed. \\
\indent The numerical results obtained in section (4) suggest that there exists a class of nonlinear problems for which the assumption $e^{t\mathcal{QL}}g(x) = g(e^{t\mathcal{QL}}x) \ \mbox{for} \ g(x) = \mathcal{QL}x$ is worth considering. 

\subsection{Computing the kernel}
The memory term for the j-th resolved variable is given by:
\begin{equation*}
\int_0^t e^{(t-s)\mathcal{L}}\mathcal{PL}F_j(x,s)ds = \int_0^t e^{(t-s)\mathcal{L}} K_j(\hat{x},s) ds = \int_0^t K_j(\hat{\phi}(t-s),s) ds
= \int_0^t K_j(\hat{\phi}(s),t-s) ds.
\end{equation*}
Both $\hat{\phi}(s)$ and $t-s$ are known quantities and therefore the challenge lies in computing  
\begin{equation*}
K_j(\hat{x},s) = \mathcal{PL} F_j(x,s).
\end{equation*}
Without the application of $\mathcal{P}$, one has:
\begin{equation}\label{eq:LF}
\mathcal{L}F_j(x,s) = \sum_{i = 1}^N R_i(x)\frac{\partial}{\partial x_i}F_j(x,s).
\end{equation}
These N sensitivities could be separately estimated, but it is more efficient to view $\mathcal{L}F$ as one single sensitivity along the unitary direction $\bar{R}(x) = R(x)/ \lVert R(x) \rVert$:
\begin{equation*}
\mathcal{L}F(x,s) = \lVert R(x) \rVert \sum_{i = 1}^N \bar{R}_i(x)\frac{\partial}{\partial x_i}F(x,s) = \lVert R(x) \rVert  ( \nabla_{\bar{R}(x)} F(x,s)) = \lVert R(x) \rVert \lim_{\epsilon -> 0} \frac{F(x + \epsilon \bar{R}(x),s) - F(x,s)}{\epsilon}.
\end{equation*}
Applying $\mathcal{P}$  and using finite difference,
\begin{equation*}
K_j(\hat{x},s) \approx  \lVert R(\hat{x}) \rVert \frac{F_j(\hat{x} + \epsilon \bar{R}(\hat{x}),s) - F_j(\hat{x},s)}{\epsilon} = \lVert R(\hat{x}) \rVert \frac{F_j(\hat{x} + \epsilon \bar{R}(\hat{x}),s)}{\epsilon}.
\end{equation*}
The memory kernel is therefore approximated as:
\begin{equation}\label{eq:shortLF}
K_j(\hat{\phi}(s),t - s) \approx  \lVert R(\hat{\phi}(s)) \rVert \frac{F_j(\hat{\phi}(s) + \epsilon \bar{R}(\hat{\phi}(s)),t - s)}{\epsilon}, \ s \in [0, \ t].
\end{equation}
In other words, evaluating the memory kernel amounts to computing the sensitivity with respect to the initial conditions, of the solution to the orthogonal dynamics in a direction that is determined by the resolved component $\hat{\phi}(s)$. Note that eq. (\ref{eq:shortLF}), and the latter statement, are independent of how the orthogonal dynamics has been solved. \\  
\indent Once the kernel is computed, the memory follows. With a rectangle rule for instance, one can make the following approximation: 
\begin{equation*}
\int_0^t K_j(\hat{\phi}(s),t-s) ds \approx \sum_{n=1}^{N_t} K_j(\hat{\phi}(s_n),t-s_n) \Delta s.
\end{equation*}

\subsection{Summary and Cost}
Denote $T$ the time domain discretized into $N_t$ equally spaced time instants $T = \{0, t_1, \ \dots, \ t_{N_t}\}$. A full-order simulation of eq. (\ref{eq:genODE}) has been run, so that the "exact" ROM solution snapshots $\{ \hat{\phi}(t_1), \ \dots \, \hat{\phi}(t_{N_t})\}$ are available. Denote $C_{full}$ the cost associated with computing these $N_t$ vector values. \\
\indent Define:
\begin{equation*}
M_n = \int_0^{t_n} K_j(\hat{\phi}(s),t-s) ds,
\end{equation*}
the memory term at $t = t_n$, which is also the subgrid term at the same instant. Then the first five terms starting at $n = 1$ can be approximated with a rectangle rule in the form:
\begin{align*}
M_1 =& \ \Delta t \ [K_j(\hat{\phi}(t_1),0)] \\
M_2 =& \ \Delta t \ [K_j(\hat{\phi}(t_1),\Delta t) \ \ + K_j(\hat{\phi}(t_2),0) ] \\
M_3 =& \ \Delta t \ [K_j(\hat{\phi}(t_1),2\Delta t) + K_j(\hat{\phi}(t_2),\Delta t) \ \ + K_j(\hat{\phi}(t_3),0) ] \\
M_4 =& \ \Delta t \ [K_j(\hat{\phi}(t_1),3\Delta t) + K_j(\hat{\phi}(t_2),2\Delta t) + K_j(\hat{\phi}(t_3),\Delta t) \ \ +  K_j(\hat{\phi}(t_4),0)] \\
M_5 =& \ \Delta t \ [K_j(\hat{\phi}(t_1),4\Delta t) + K_j(\hat{\phi}(t_2),3\Delta t) + K_j(\hat{\phi}(t_3),2\Delta t) +  K_j(\hat{\phi}(t_4),\Delta t) + K_j(\hat{\phi}(t_5),0)].
\end{align*}
Thus, to compute  $M_1$ up to $M_5$, one must compute:
\begin{itemize}
    \item $K_j(\hat{\phi}(t_1),s)$ for $s \in \{0,\ \Delta t,\ 2\Delta t,\ 3\Delta t,\ 4\Delta t\}$;
    \item $K_j(\hat{\phi}(t_2),s)$ for $s \in \{0,\ \Delta t,\ 2\Delta t,\ 3\Delta t\}$;
    \item $K_j(\hat{\phi}(t_3),s)$ for $s \in \{0,\ \Delta t,\ 2\Delta t\}$;
    \item $K_j(\hat{\phi}(t_4),s)$ for $s \in \{0,\ \Delta t\}$;
    \item $K_j(\hat{\phi}(t_4),s)$ for $s \in \{0 \}$;
\end{itemize}

According to eq. (\ref{eq:shortLF}), each kernel value $K_j(\hat{\phi}(t_n),s)$ requires computing $F_j(\hat{\phi}(t_n) + \epsilon \bar{R}(\hat{\phi}(t_n)),s)$. It is obtained - see section 3.1. - by solving the pseudo orthogonal ODE:
\begin{equation}\label{eq:spawnOrtODE}
\frac{\partial}{\partial t} \phi^{\mathcal{Q}}(t) = R(\phi^{\mathcal{Q}}(t)) - R(\hat{\phi}^{\mathcal{Q}}(t))\ \ ; \ \ \phi^{\mathcal{Q}}(0) = \hat{\phi}(t_n) + \epsilon \bar{R}(\hat{\phi}(t_n)),
\end{equation}
and retrieving the first $m$ components of the right-hand side term at $t = s$. \\
\indent Let's assume that the costs of evolving eq. (\ref{eq:spawnOrtODE}) and eq. (\ref{eq:genODE}) in time are the same. Then the cost $C_{M-Z}$ of computing the memory terms $M_n$ at all $N_t$ time instants scales as:
\begin{equation}\label{cmzeqn}
C_{M-Z} \ \propto \ \frac{N_t-1}{2}C_{full}.
\end{equation}
\indent Since computing $M_n$ requires computing $n$ values of $K_j(\hat{\phi}(s), t_n - s)$ the scaling that is obtained should not be surprising. It is emphasized again that the purpose of the proposed method is to estimate the memory kernel a priori in nonlinear systems.  
\\
\indent Such scaling arguably limits the size of the problems that can be considered. Note however, that the described procedure can be easily implemented in parallel. The kernel values $K_j(\hat{\phi}(t_n),s)$ at the required instants $s$ can be computed separately for each $n = 1 \dots N_t $. \\
\indent Eq. (\ref{cmzeqn}) gives the cost of computing the kernel values in the full range $s \in [0, \ t]$. If the size of the problem is such that $C_{MZ}$ is a limiting factor, one can still use the procedure to compute the kernel $K_j(\hat{\phi}(s), t - s)$ in a truncated range $s \in [t - \tau, \ t]$ with $\tau$ taken such that the corresponding cost $C_{MZ}(\tau)$ is affordable. Assuming that the method is always accurate, this could be seen as an indirect way of investigating finite memory. 
For all the problems considered in this work, the kernel was estimated in the full range  $s \in [0, \ t]$ even though the observations that followed suggested that it was not needed. 

\section{Applications}
In this section, the effectiveness of the kernel estimation technique is assessed in three different problems - the Brusselator, the Viscous Burgers equation and the Kuramoto-Sivashinsky equation. 
\subsection{Brusselator}

\begin{figure}[h]
    \centering
    \subfigure[(A,B) = (1, 1.7)]{\includegraphics[width = .48\linewidth, height = 0.30\textheight]{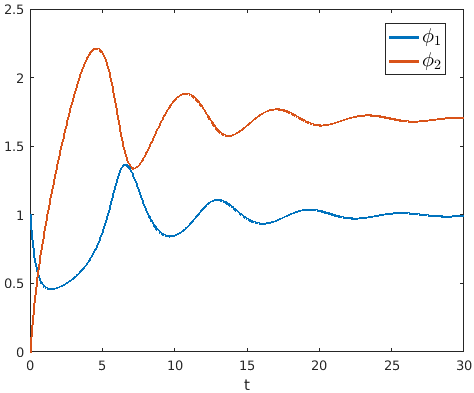}}
    \subfigure[(A,B) = (1, 3)]{\includegraphics[width = .48\linewidth, height = 0.30\textheight]{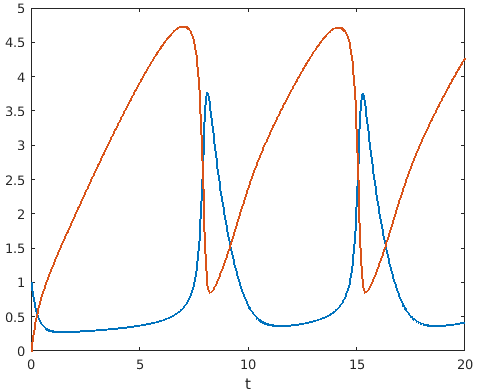}}\\
    \caption{Brusselator: Full order solutions for two distinct configurations and $(\phi_1, \phi_2) = (1, 0)$.}
    \label{fig:intro_bru}
\end{figure}

\indent Two configurations of the Brusselator eq. (\ref{eq:Bru})  are considered (see figure \ref{fig:intro_bru}):
\begin{itemize}
    \item $(A,B) = (1,1.7)$, where the system is asymptotically stable.
    \item $(A,B) = (1,3)$, where the system exhibits Limit Cycle Oscillations (LCO).
\end{itemize}
\indent The Liouville operator for the Brusselator is given by:
\begin{equation*}
\mathcal{L} \equiv (A + x_1^2 x_2 - B x_1 - x_1)\frac{\partial}{\partial x_1} + (B x_1 - x_1^2 x_2)\frac{\partial}{\partial x_2}.
\end{equation*}

\noindent The pseudo orthogonal ODE is given by:
\begin{equation}\label{eq:BruOrtODE}
\begin{cases}
\begin{gathered}
\frac{d}{dt} \phi_1^{\mathcal{Q}} = \ \ \ (\phi_1^{\mathcal{Q}})^2 \phi_2^{\mathcal{Q}}, \\
\frac{d}{dt} \phi_2^{\mathcal{Q}} = -(\phi_1^{\mathcal{Q}})^2 \phi_2^{\mathcal{Q}}, \\
\phi^{\mathcal{Q}}(0) = x,
\end{gathered}
\end{cases}
\end{equation}
and $F_1(x,s) = \phi_1^{\mathcal{Q}}(x,s)^2 \phi_2^{\mathcal{Q}}(x,s)$. \\ 
\indent The original ODE (\ref{eq:Bru}) and the pseudo orthogonal ODE (\ref{eq:BruOrtODE}) are solved using $4^{th}$ order Runge-Kutta time integration.

\begin{figure}[h]
    \centering
    \subfigure[(A,B) = (1, 1.7) - $\Delta t = 1.5\ 10^{-2}$ s]{\includegraphics[width = .48\linewidth, height = 0.30\textheight]{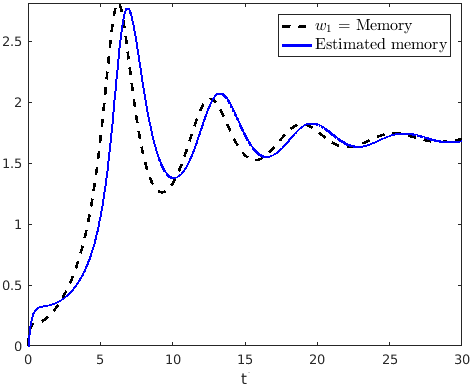}}
    \label{fig:Bru_M-Z_stable}
    \subfigure[(A,B) = (1, 3) - $\Delta t = 10^{-2} $s]{\includegraphics[width = .48\linewidth, height = 0.30\textheight]{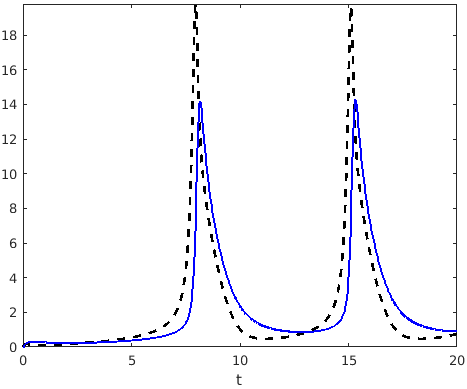}}
    \label{fig:Bru_M-Z_LCO}
    \caption{Brusselator: The dotted curve is $w_1(t) = \phi_1(t)^2 \phi_2(t)$ in eq. (\ref{eq:MZref}) computed using full-order solution snapshots. The plain curve is the same quantity but evaluated by integrating the kernel (RHS of eq. (\ref{eq:MZref}))  over time. 
    }
    \label{fig:BRUMZ}
\end{figure}

\begin{figure}[h]
    \centering
    \subfigure[t = 4.28 s]{\includegraphics[width = .48\linewidth, height = 0.24\textheight]{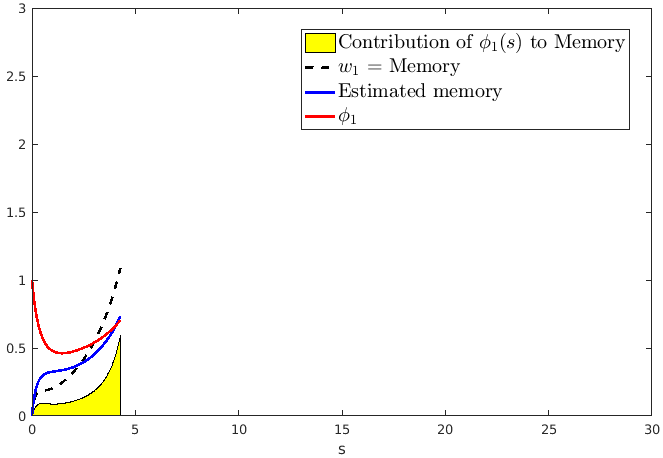}}
    \subfigure[t = 8.55 s]{\includegraphics[width = .48\linewidth, height = 0.24\textheight]{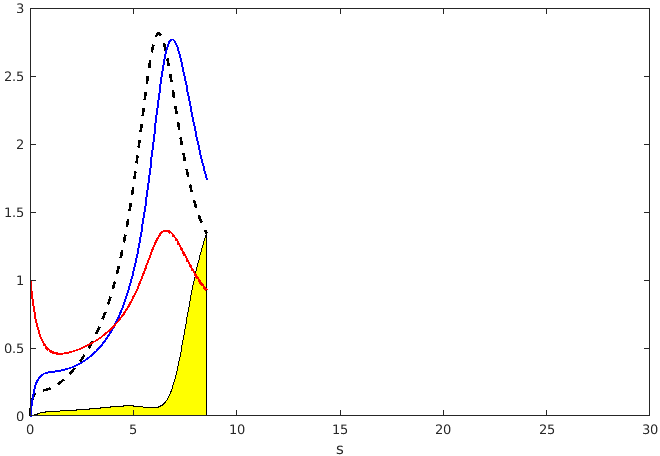}}
    \subfigure[t = 12.83 s]{\includegraphics[width = .48\linewidth, height = 0.24\textheight]{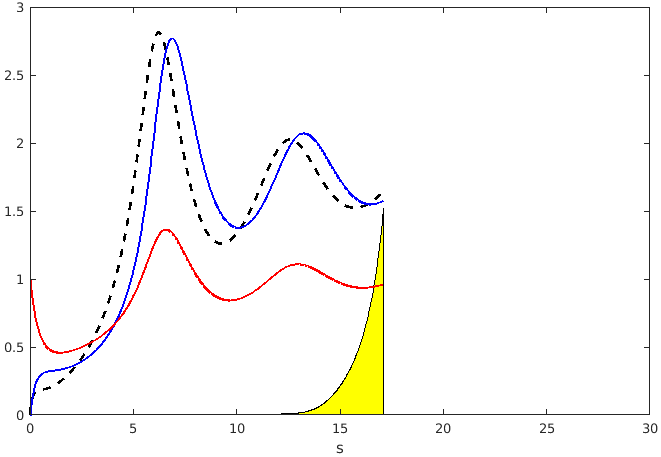}}
    \subfigure[t = 30 s]{\includegraphics[width = .48\linewidth, height = 0.24\textheight]{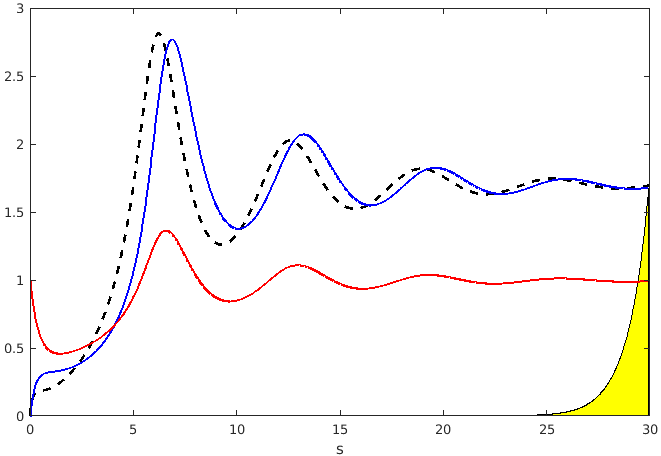}}
    \caption{Brusselator (stable configuration): Memory kernel snapshots. $w_1(t) = \phi_1(t)^2 \phi_2(t)$. The (estimated) contribution of $\phi_1(s)$ to the memory at $t$ (yellow shaded area) is $K_1(\phi_1(s),t-s) ds$ for $s \in [0, \ t]$.}
   \label{fig:Bru_stable_decay}
\end{figure}

\begin{figure}[h]
   \centering
    \subfigure[t = 8 s] {\includegraphics[width = .48\linewidth, height = 0.24\textheight]{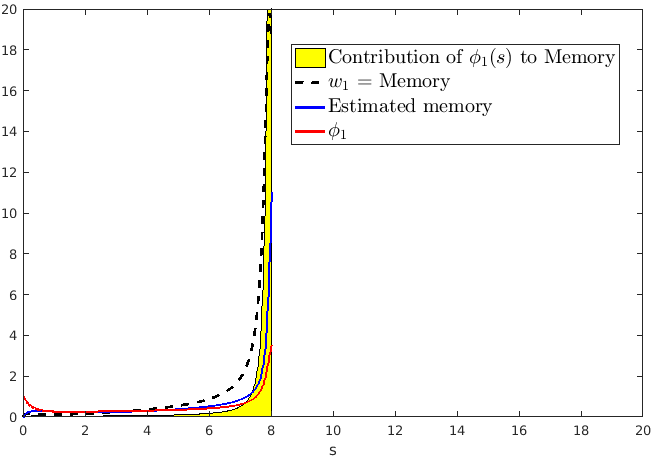}}
    \subfigure[t = 9 s] {\includegraphics[width = .48\linewidth, height = 0.24\textheight]{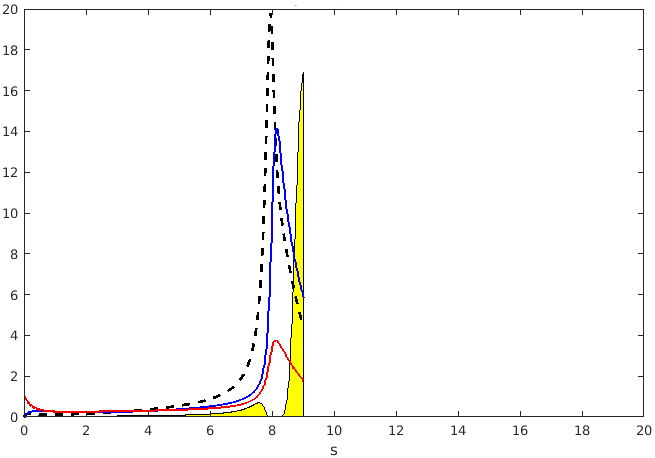}}
    \subfigure[t = 10 s]{\includegraphics[width = .48\linewidth, height = 0.24\textheight]{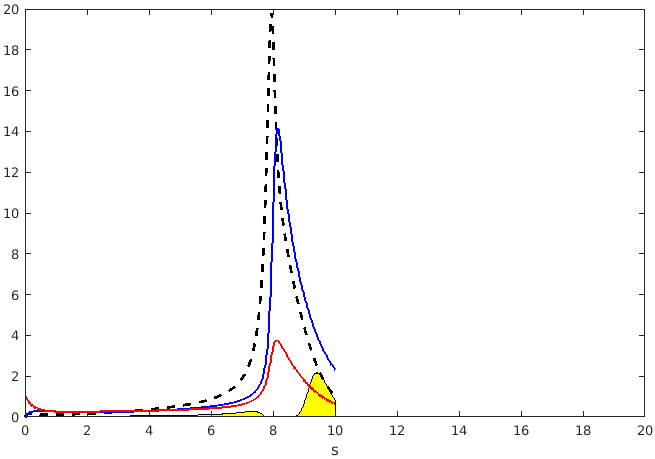}}
    \subfigure[t = 15 s]{\includegraphics[width = .48\linewidth, height = 0.24\textheight]{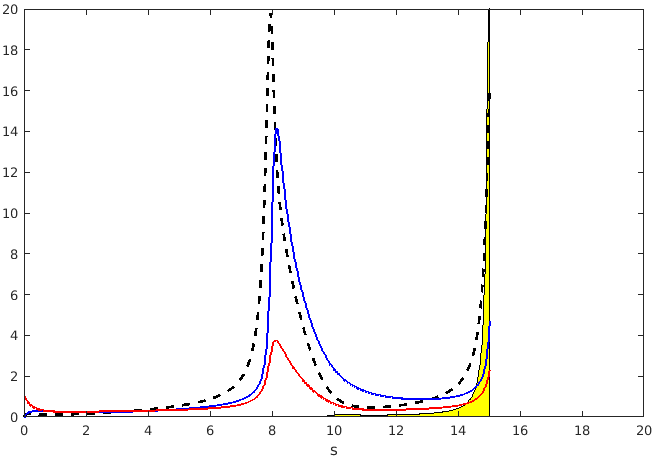}}
    \caption{Brusselator (LCO configuration): Memory kernel snapshots.}
    \label{fig:Bru_decay}
\end{figure}

The kernel $K_1$ is then computed at the required points using eq. (\ref{eq:shortLF}) 
with $\epsilon = 10^{-8}$. Figures \ref{fig:BRUMZ}.(a) and \ref{fig:BRUMZ}.(b) show the numerical results for the stable and LCO configurations, respectively. There are some amplification and phase errors, but the results are, in general, good. \\
\indent In figures \ref{fig:Bru_stable_decay} and \ref{fig:Bru_decay} we provide snapshots of the estimated memory kernel $K_1(\phi_1(s), t - s)$ (the same kind of data as in section (2.1.)) for $(x_1,x_2) = (1,0)$. Figures \ref{fig:Bru_stable_decay}(a)-(d) were computed for the stable configuration $(A,B) = (1,1.7)$ while figures \ref{fig:Bru_decay}(a)-(d) were computed for the unstable configuration $(A,B) = (1,3)$. In these figures, we see that, in both configurations, the memory has finite support. 

\subsection{Viscous Burgers Equation}
The Viscous Burgers equation in discrete Fourier space is considered:
\begin{equation}\label{eq:Burgers}
\frac{d}{dt} u_k + \frac{ik}{2}\sum_{\substack{p+q = k \\ (p,q) \in F \cup G}} u_p u_q = -\nu k^2 u_k, \ k \in F \cup G.
\end{equation}

\indent The physical domain $[0, \ 2\pi]$ is discretized into $2 N$ equidistant points $\{ x_0, \ x_1, \ \dots , x_{2 N-1}\}$ with $x_0 = 0$ and $x_{2 N-1} = 2 \pi - \Delta x$. The parameter $\nu$ is the viscosity.
The relationship between the coefficients $u_k$ and the grid values $u(x_n)$ is given by: 
\begin{equation*}
u_k = \sum_{n = 0}^{2N - 1} u(x_n) e^{2\pi i k \frac{n}{N}}.
\end{equation*}
F and G represent the sets of resolved and unresolved Fourier modes, respectively. Denote $\bm{u} = (u_k)_{1 \leq k \leq N}$ the vector of Fourier coefficients and $\bm{\hat{u}}$ the same vector with all of G modes set to 0. Denote $\bm{u}_0 = \bm{u}(0)$. The initial conditions are given by \cite{Wang}:
\begin{equation*}
u(x,0) = \   \sum_{k = 1}^{m} \sqrt[]{2E(k)}\sin(k x + \beta_k)\ \ ; \  E(k) = \left\{
    \begin{array}{ll}
        5^{-5/3} & \mbox{if } 1 \leq k \leq 5 \\
        k^{-5/3} & \mbox{if } k > 5
    \end{array}
\right.
\end{equation*}
The subgrid terms for this problem are given by:
\begin{equation*}
w_k = - \frac{ik}{2}\big[\sum_{\substack{p+q = k \\ (p,q) \in F \cup G}} u_{p} u_{q} - \sum_{\substack{p+q = k \\ (p,q) \in F^2}} u_{p} u_{q} \ \Big ] = \mathcal{P}\int_0^t K_k(\bm{\hat{u}}(s),t-s)ds, \ k \in F.
\end{equation*}
The Liouville operator is:
\begin{equation*}
\mathcal{L} \equiv \sum_{k = 1}^{N-1} (-\frac{ik}{2} \sum_{\substack{p+q = k \\ (p,q) \in F \cup G}} u_{0p}u_{0q} - \nu k^2 u_{0k})\frac{\partial}{\partial u_{0k}}.
\end{equation*}
Applying eq.~(\ref{eq:shortLF}) to the VBE, we obtain:
\begin{equation*}
K_j(\bm{\hat{u}}(s),t-s) \approx \lVert R(\bm{\hat{u}}(s)) \rVert \frac{F_j(\bm{\hat{u}}(s) + \epsilon \bar{R}(\bm{\hat{u}}(s)),t-s)}{\epsilon}.
\end{equation*}
Where $F(\bm{u_0},s)$ is the RHS term at time $s$ of the corresponding pseudo orthogonal ODE:
\begin{equation}
\label{eq:ortODEBurgers}
\begin{cases}
\begin{gathered}
\frac{d}{dt} u_k^{\mathcal{Q}} = - \frac{ik}{2}[\sum_{\substack{p+q = k \\ (p,q) \in F \cup G}} u_{p}^{\mathcal{Q}} u_{q}^{\mathcal{Q}} - \sum_{\substack{p+q = k \\ (p,q) \in F^2}} u_{p}^{\mathcal{Q}} u_{q}^{\mathcal{Q}}], \hspace{2.4 cm} k \in F \\
\frac{d}{dt} u_k^{\mathcal{Q}} = - \frac{ik}{2}[\sum_{\substack{p+q = k \\ (p,q) \in F \cup G}} u_{p}^{\mathcal{Q}} u_{q}^{\mathcal{Q}} - \sum_{\substack{p+q = k \\ (p,q) \in F^2}} u_{p}^{\mathcal{Q}} u_{q}^{\mathcal{Q}}] - \nu k^2 u_k^{\mathcal{Q}}, \hspace{1 cm} k \in G \\
\bm{u^{\mathcal{Q}}}(0) = \bm{u_0}.
\end{gathered}
\end{cases}
\end{equation}
\subsubsection{Results}
\begin{figure}[!htbp]
    \centering
     \subfigure[Integration over time of the estimated kernel]{\includegraphics[width = .49\linewidth, height = 0.36\textheight]{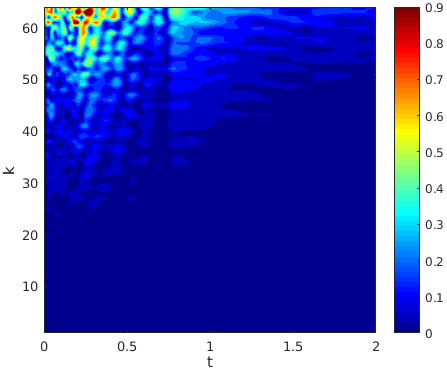}}
     \subfigure[Using the full-order solution]{\includegraphics[width = .49\linewidth, height = 0.36\textheight]{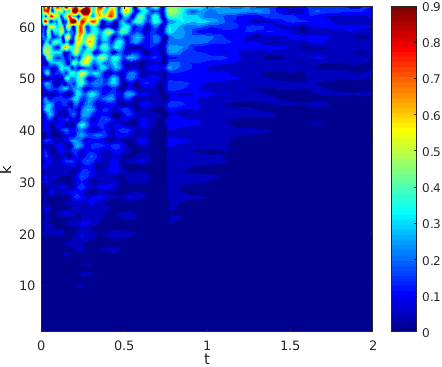}}
    \caption{Viscous Burgers equation: Norm of $(w_{k})_{1 \leq k \leq 64}$ represented as a continuous contour.}    
    \label{fig:BMZresults_contour}
\end{figure}

\begin{figure}[!htbp]
    \centering
    \subfigure[mean $(w_{k})_{1 \leq k \leq Nc}$ - real part]{\includegraphics[width = .47\linewidth, height = 0.36\textheight]{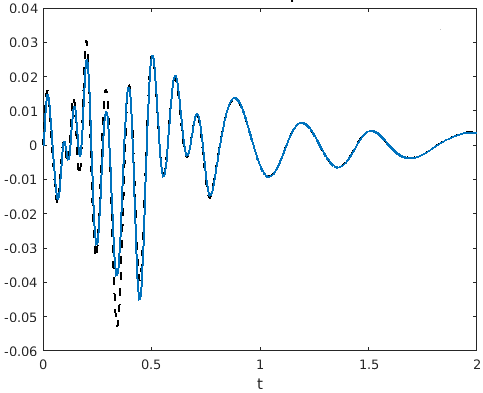}}
    \subfigure[Contribution of the subgrid terms to the resolved energy decay]{\includegraphics[width = .47\linewidth, height = 0.36\textheight]{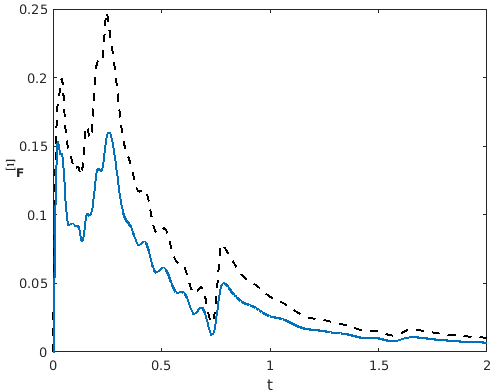}}
    \caption{Viscous Burgers equation: Numerical results (dotted - using the full-order solution; plain - integration of the estimated kernel) under the two metrics.}
    \label{fig:BMZresults}
\end{figure}

\begin{figure}[!h]
\centering
    \subfigure[$w_{61}$]{\includegraphics[width = .48\linewidth, height = 0.24\textheight]{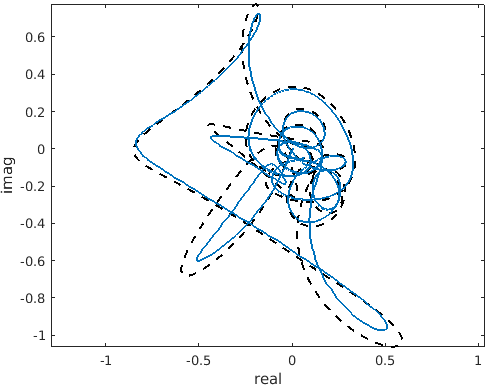}}
    \subfigure[$w_{62}$]{\includegraphics[width = .48\linewidth, height = 0.24\textheight]{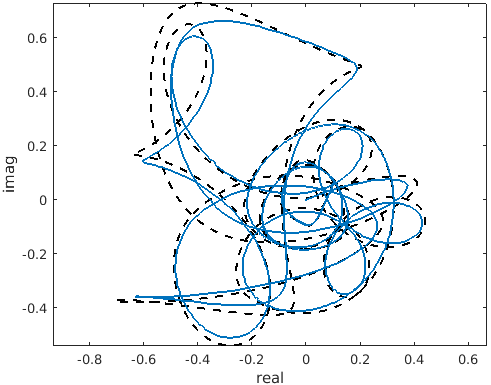}}
    \subfigure[$w_{63}$]{\includegraphics[width = .48\linewidth, height = 0.24\textheight]{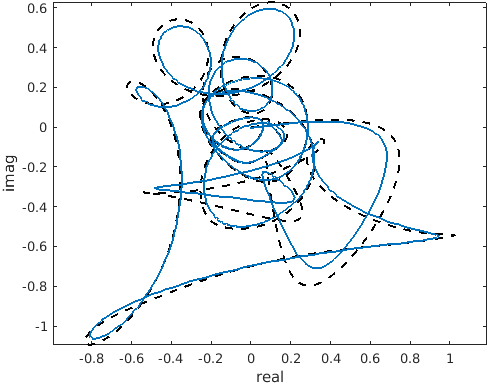}}
    \subfigure[$w_{64}$]{\includegraphics[width = .48\linewidth, height = 0.24\textheight]{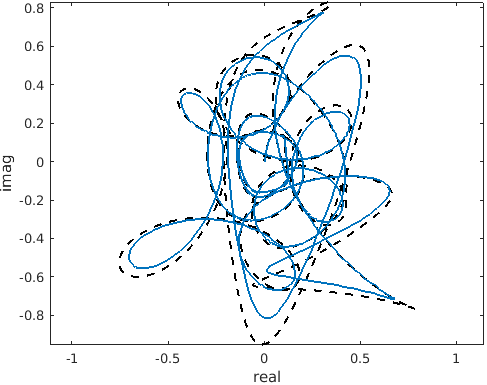}}
    \caption{Viscous Burgers: subgrid/memory terms (dotted - using the full-order solution snapshots, plain - integrating over time the estimated kernel).}
    \label{fig:BModes}
\end{figure}

\begin{figure}[!htbp]
\centering
    \subfigure[t = 0.5]{\includegraphics[width = .46\linewidth, height = 0.22\textheight]{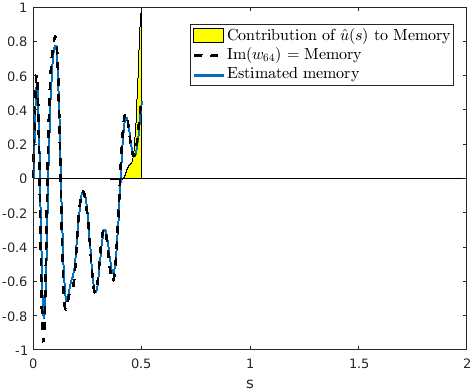}}
    \subfigure[t = 1.0]{\includegraphics[width = .46\linewidth, height = 0.22\textheight]{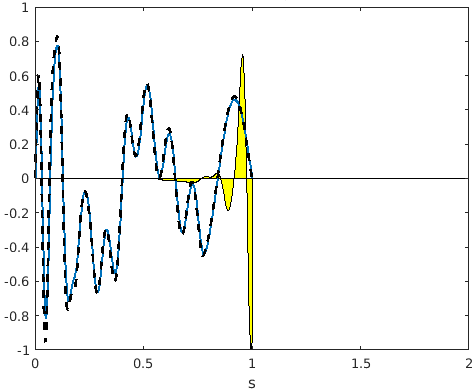}}
    \subfigure[t = 1.5]{\includegraphics[width = .46\linewidth, height =0.22\textheight]{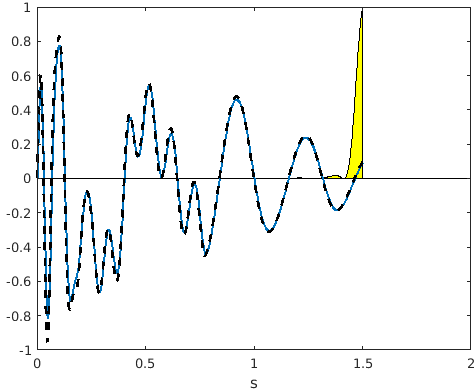}}
    \subfigure[t = 2.0]{\includegraphics[width = .46\linewidth, height = 0.22\textheight]{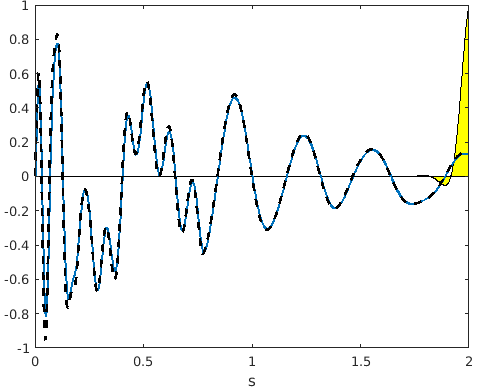}}
    \caption{Viscous Burgers. Snapshots of the kernel (imaginary part) for cut-off mode 64. }
    \label{fig:BSnaps32}
\end{figure}

\begin{figure}[!htbp]
\centering
    \subfigure[decay profiles]{\includegraphics[width = .48\linewidth, height = 0.32\textheight]{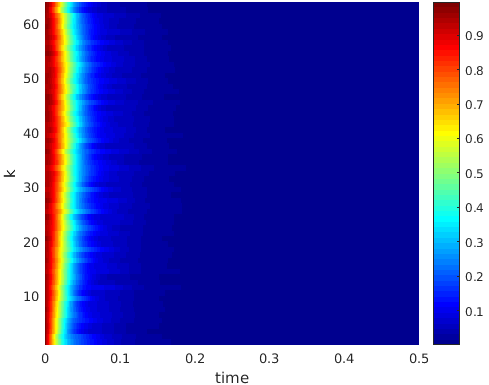}}
    \subfigure[Corresponding memory lengths (99\%) ]{\includegraphics[width = .48\linewidth, height = 0.32\textheight]{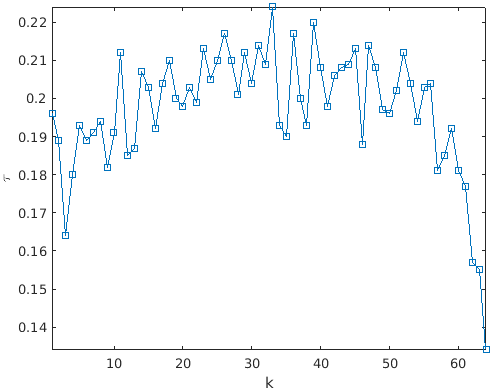}}
    
    \caption{Viscous Burgers. Figure (a) provides a first estimation of how fast the kernel $K_j(\hat{u},t)$ can be expected to decay with time $t$. Figure (b) provides, for each wavenumber k, an estimation of the time $\tau$ after which the decay profile goes below 1\%.}
    \label{fig:Burgdecay}
\end{figure}

\begin{figure}[!htbp]
\centering
    \subfigure['cut-off' decay profiles for various values of $m$]{\includegraphics[width = .48\linewidth, height = 0.35\textheight]{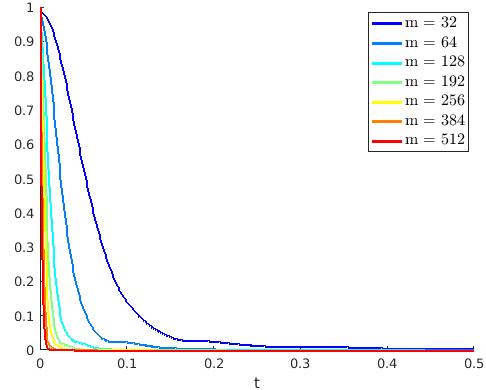}}
    \subfigure[Corresponding memory lengths (99\%)]{\includegraphics[width = .48\linewidth, height = 0.35\textheight]{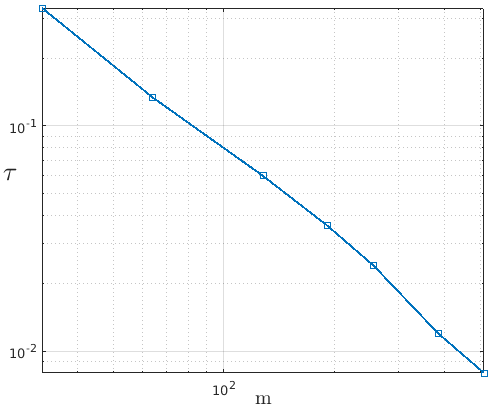}}
    \caption{Viscous Burgers: First hindsight into the effects of the ROM size $m$ on the memory length. The bigger the ROM faster the kernel $K_{m}$ decays. Figure (a) shows averaged cut-off decay profiles $K_{m}(\hat{u},t)$ for different values of $m$. Figure (b) hints at a power-law relationship between $m$ and the memory length.}
    \label{fig:BdecayNs}
\end{figure}

The parameters are $m = 64$, $N = 1048$, $t_f = 2$, $\epsilon = 10^{-8}$, $dt = 10^{-3}s$ and $\nu = 10^{-3}$. A fourth-order Runge-Kutta integration scheme with 3/2 padding is used to solve  eqns. (\ref{eq:Burgers}) and (\ref{eq:ortODEBurgers}). \\
\indent A global picture of the results is given in figure \ref{fig:BMZresults_contour} where contour plots of the norm of the subgrid terms $(w_k)_{1\leq k \leq 64}$ computed using the full-order solution and computed with our method are given. Two metrics are used to measure how well the subgrid terms were reproduced: 

\begin{itemize}
\item The first metric is a mean over all the resolved modes (figure \ref{fig:BMZresults}(a)). The agreement between the estimation and the full-order results is excellent, although amplification errors are noticeable in the first third of the modes. This metric only quantifies whether the dominant memory terms have been well reproduced. For this fundamental problem in particular, the higher the wavenumber $k$, the higher $w_k$. 

\item To obtain a more complete picture, the second metric analyzes the contribution of the subgrid terms to the resolved part of the energy decay rate. Since 
$$\left[-\frac{dE}{dt}\right]_F=-\frac{1}{2}\sum_{k \in F} (u_k \frac{d}{dt}u_k^{*} + u_k^*\frac{d}{dt}u_k) = - \sum_{k \in F} Re(u_k^* \frac{d}{dt}u_k),$$

the contribution from $w_k$ is defined as $\Xi_F = - \sum_{k \in F} Re(u_k^* w_k)$. The fact that the low wavenumber fourier modes $u_k$ are more energetic than the high wavenumber modes compensates for the opposite trend in the subgrid terms $w_k$, and makes $\Xi_F$ a good indicator of whether the low wavenumber subgrid terms have been well reproduced (figure \ref{fig:BMZresults} (b)). It appears that for this problem, the proposed method only enables a faithful reproduction of the high wavenumber subgrid terms. 
\end{itemize}
\indent Figure \ref{fig:BModes} features the last four subgrid/memory terms. Figure \ref{fig:BSnaps32} shows snapshots of the imaginary part of the memory kernel $K_{64}(\hat{u}(s),t-s)$ at four different time instants $t$. These figures reveal the finite support of the memory length and a 
decaying sinusoidal form of the kernel. \\
\indent To estimate the memory length, the profiles for the kernels $K_j$ have been built by averaging and scaling $|K_j(\hat{\phi}(t_n),s)|$ over all the trajectory points. This is shown in figure \ref{fig:Burgdecay}.(a) as a discrete contour plot. Figure \ref{fig:Burgdecay}.(b) gathers the memory length for those $m$ decay profiles. 
Figures~\ref{fig:BdecayNs}(a-b) show the effect of the size of the coarse-grained model the estimated memory length for $w_{64}$ ($k = 64$ is the cut-off wavenumber). As expected, it is seen that, as the number of resolved modes increases,  the memory length diminishes, as expected. For this specific problem, the relationship between the memory length and the level of coarse-graining amounts to a simple scaling law.

\subsection{Kuramoto-Sivashinsky (K-S) equation}
\begin{figure}
\centering
\includegraphics[scale = 0.85]{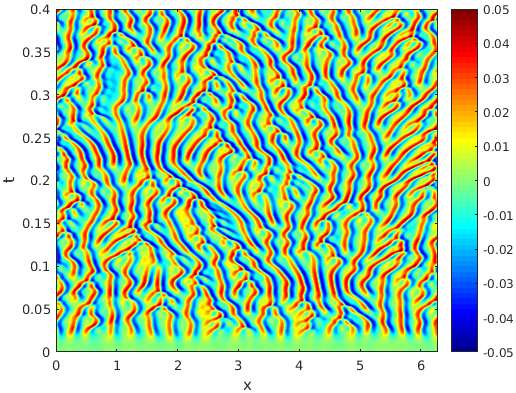}
\caption{Kuramoto-Sivashinsky: Full-order solution in physical space.}
\label{fig:KScontour}
\end{figure}

The K-S equation~\cite{Hyman, Nicolaenko} in discrete Fourier space is given by:
\begin{equation*}
\frac{d}{dt}u_k = (k^2 - \nu k^4) u_k - \frac{ik}{2}\sum_{\substack{p+q = k}}u_p u_q, \ k \in F \cup G.
\end{equation*}
\indent This equation has similarities with the VBE, but develops richer behaviour due to the presence of both second and fourth order spatial derivatives. The linear part $k^2 - \nu k^4$ is strictly positive for $k < 1/\sqrt[]{\nu}$. The corresponding modes act as sources of energy for the entire system. Stability is ensured by the convection term that transfers the energy produced by the source at the large scales (low wavenumber) to the small scales (high wavenumber) where $-\nu k^4$ dominates. The viscosity $\nu$ is chosen such that the system would display (figure \ref{fig:KScontour}) a state of persistent dynamical disorder. The  initial conditions are assumed to be the same as in the Viscous Burgers equation. 
Since the K-S system is stiffer than the viscous Burgers equation,  a fourth order Exponential Time Differencing (ETD) scheme~\cite{Cox,Kassam} is used instead of the Runge-Kutta scheme for robustness.  The initial conditions and problem parameters are the same as for Burgers equation, except $m = 32$, $t_f = 0.4$ s and $dt = 10^{-4}$s. \\

\subsubsection{Results}
\begin{figure}[!htbp]
    \centering
    \subfigure[Integration over time of the estimated kernel.]{\includegraphics[width = .97\linewidth, height = 0.37\textheight]{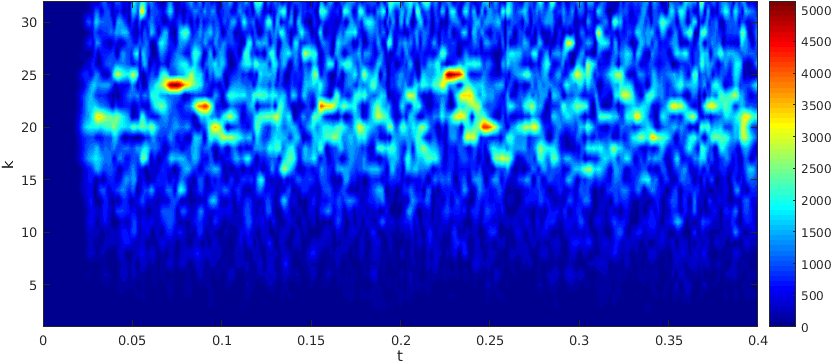}}
    \subfigure[Computed from the full-order solution.]{\includegraphics[width = .97\linewidth, height = 0.37\textheight]{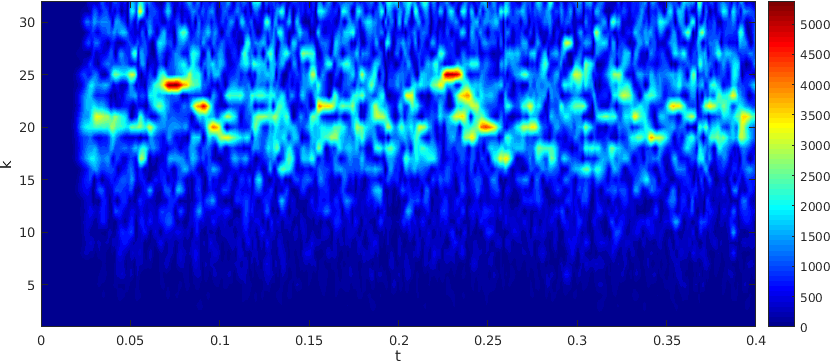}}
    \caption{Kuramoto-Sivashinsky: Norm of $(w_{k})_{1 \leq k \leq 32}$ represented as a continuous contour.}
    \label{fig:KSResults}
\end{figure}

\begin{figure}[!htbp]
    \centering
    \includegraphics[width = .97\linewidth, height = 0.37\textheight]{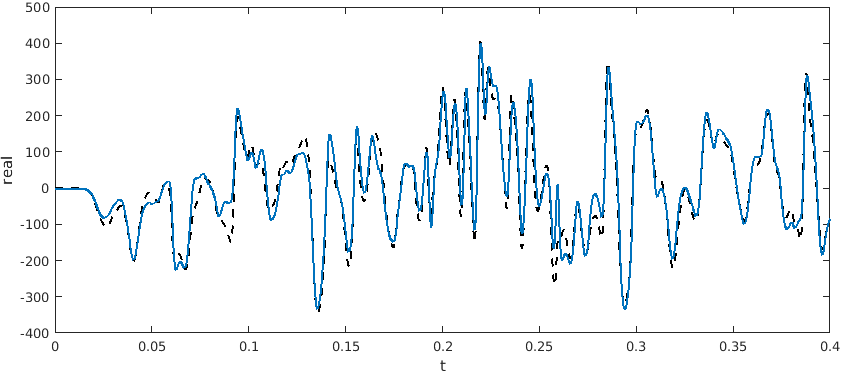}
    \caption{Kuramoto-Sivashisky. Numerical results under the first metric - mean $(w_{k})_{1 \leq k \leq Nc}$ (real part) versus time. The dotted line is computed using the full-order solution. The plain line is obtained by integrating over time the estimated kernel.}
    \label{fig:KSmeanRe}
\end{figure}

\begin{figure}[!htbp]
    \centering
    \includegraphics[width = .97\linewidth, height = 0.37\textheight]{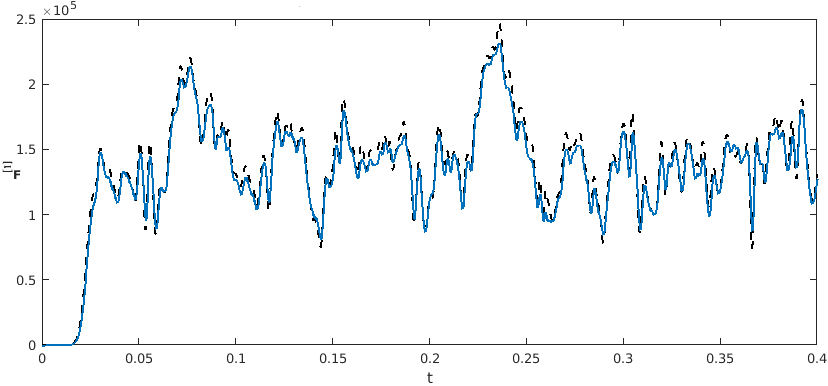}
    \caption{Kuramoto-Sivashinksy. Numerical results under the second metric $\Xi_F$. Same legend as in figure \ref{fig:KSmeanRe}.}
    \label{fig:KSXi}
\end{figure}

\begin{figure}[!h]
\centering
    \subfigure[$w_{21}$ - Exact]{\includegraphics[width = .22\linewidth, height = 0.22\textheight]{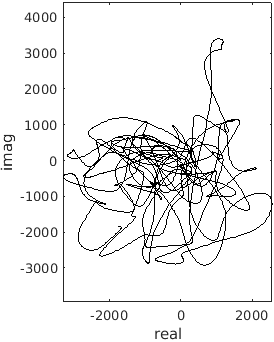}}
    \subfigure[$w_{21}$ - MZ]{\includegraphics[width = .22\linewidth, height = 0.22\textheight]{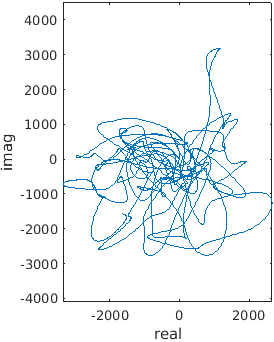}}
    \subfigure[$w_{25}$ - Exact]{\includegraphics[width = .22\linewidth, height = 0.22\textheight]{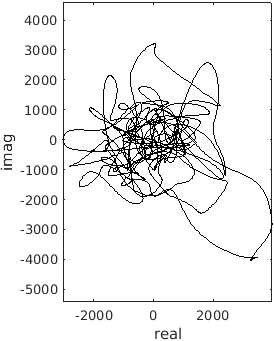}}
    \subfigure[$w_{25}$ - MZ]{\includegraphics[width = .22\linewidth, height = 0.22\textheight]{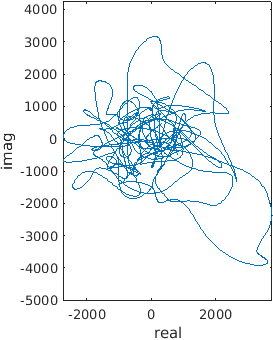}}
    \subfigure[$w_{29}$ - Exact]{\includegraphics[width = .22\linewidth, height = 0.22\textheight]{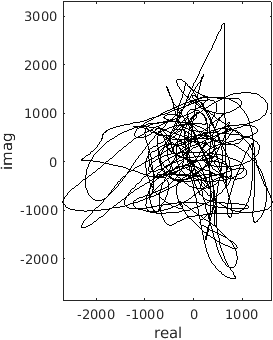}}
    \subfigure[$w_{29}$ - MZ]{\includegraphics[width = .22\linewidth, height = 0.22\textheight]{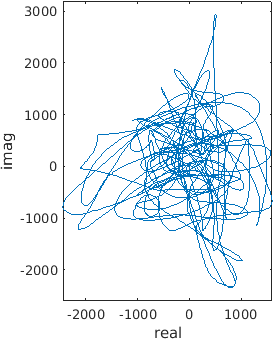}}
    \subfigure[$w_{32}$ - Exact]{\includegraphics[width = .22\linewidth, height = 0.22\textheight]{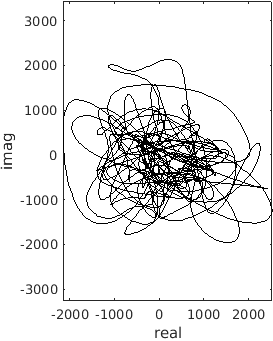}}
    \subfigure[$w_{32}$ - MZ]{\includegraphics[width = .22\linewidth, height = 0.22\textheight]{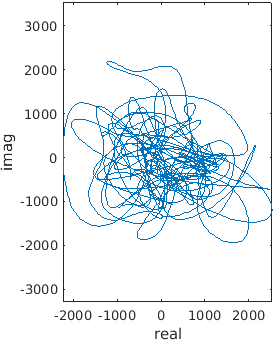}}
    \caption{Kuramoto-Sivashinsky. Subgrid/memory terms ("Exact" - using the full-order solution snapshots, "MZ" - integrating over time the estimated kernel).}
    \label{fig:KSmodes}
\end{figure}

\begin{figure}[!htbp]
\centering
    \subfigure[t = 0.1]{\includegraphics[width = .46\linewidth, height = 0.22\textheight]{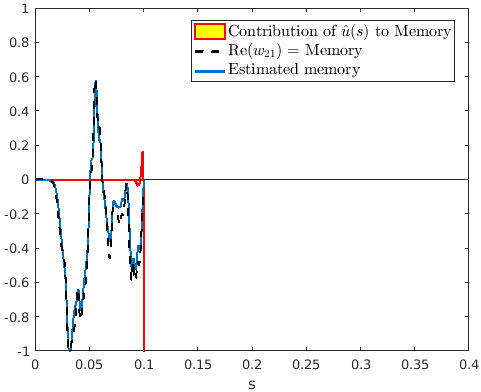}}
    \subfigure[t = 0.2]{\includegraphics[width = .46\linewidth, height = 0.22\textheight]{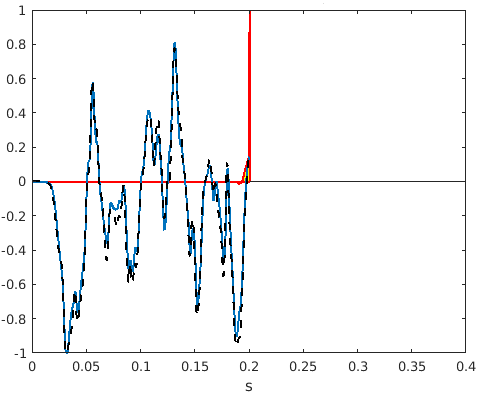}}
    \subfigure[t = 0.3]{\includegraphics[width = .46\linewidth, height = 0.22\textheight]{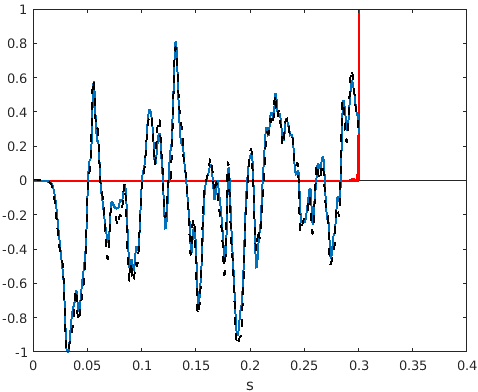}}
    \subfigure[t = 0.4]{\includegraphics[width = .46\linewidth, height = 0.22\textheight]{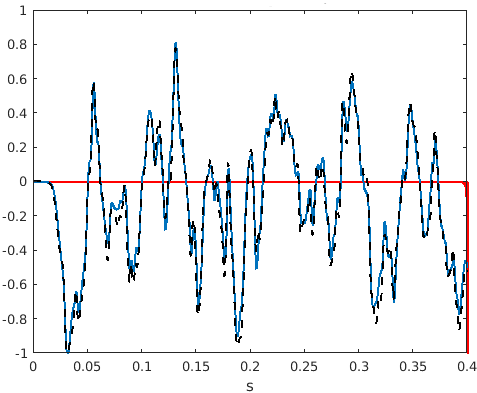}}
    \caption{Kuramoto-Sivanshinksy. Snapshots of the kernel (real part) for mode 21. The estimated kernels decay very rapidly. The memory at time $t$ is completely determined by the most recent kernel values $K_j(\hat{u}(s),t-s), \ s \in [t- \tau, \ t]$. Note: The memory terms and the computed kernels have been scaled to fit together in the graphs. }
    \label{fig:KS31}
\end{figure}

\begin{figure}[!htbp]
\centering
    \subfigure[Decay profiles]{\includegraphics[width = .48\linewidth, height = 0.30\textheight]{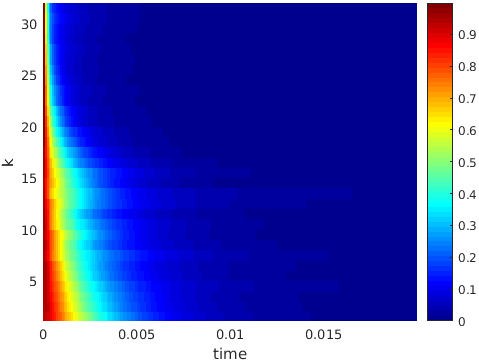}}
    \subfigure[Corresponding memory lengths (99\%) ]{\includegraphics[width = .48\linewidth, height = 0.30\textheight]{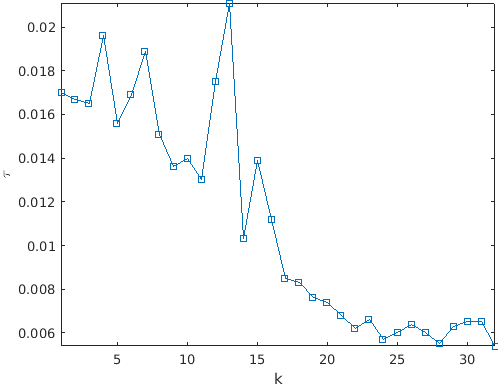}}
    
    \caption{Kuramoto-Sivashinsky. Figure (a) provides a first estimation of the rate at which the kernel $K_j(\hat{u},s)$  decays with time. Figure (b) provides, for each wavenumber k, an estimation of the time $\tau$ after which the decay profile goes below 1\% of its maximal amplitude.}
    \label{fig:KSdecay}
\end{figure}

The overall picture (contour of the subgrid terms) is given in figure \ref{fig:KSResults}. The mean of $w_k$ and $\Xi_F$ are shown in figures \ref{fig:KSmeanRe} and \ref{fig:KSXi}, respectively. By both metrics, the method performs well. Figure \ref{fig:KSmodes} compares the subgrid terms $w_j(t)$ computed using the full-order solutions and the estimated memory terms for selected modes $j$. Figure \ref{fig:KS31} suggests a very short memory length compared to the Burgers equation. The estimated decay profiles and their corresponding memory lengths are provided in figures \ref{fig:KSdecay}.(a) and \ref{fig:KSdecay}.(b).\\
\indent The presence of a finite support in time for the memory suggests that, in order to reconstruct the subgrid terms $w_j$ at time $t$ from the M-Z perspective, one would only need to compute the kernel $K_j(\hat{\phi}(s), t - s)$ for $s \in [t - \tau, \ t]$.


\section{Conclusion and perspectives}
The Mori-Zwanzig (M-Z) formalism provides a mathematically consistent framework for the construction of coarse-grained models of dynamical systems. Using a projection operator that partitions the state vector into resolved and unresolved parts, the M-Z approach recasts a high-order nonlinear ODE system into Generalized Langevin Equation (GLE) for the resolved part. The  unresolved part is exactly represented in the form of a non-local memory integral that involves the past history of the resolved physics. 

In this work, we presented an a priori procedure to estimate the memory kernel for the truncation projector. In this context, computing the memory kernel $K_j(\hat{\phi}(s),t-s)$ amounts to the evaluation of the sensitivity with respect to the initial conditions of the solution to the orthogonal dynamics a direction imposed by the ROM trajectory. 

Instead of solving the orthogonal dynamics equation - which is a high-dimensional linear partial differential equation - we propose  to work with a simplified equation. This simplification is a result of the assumption that the semi-group of the orthogonal dynamics is a composition operator for the observable $g(x) = \mathcal{QL}x$. The ersatz has the advantage of being tractable by evolving in time a characteristic equation that we call the pseudo orthogonal ODE. It has to be emphasized that the method and analysis proposed in this work is limited to the truncation projector only.\\
\indent The procedure is exact for linear systems for which the kernel is known analytically. Based on numerical results for the Brusselator, the Viscous Burgers equation and the K-S equation, we conjecture that the assumption made on the semi-group of the orthogonal dynamics is appropriate for a certain class of nonlinear problems. For Burgers, the dominating, near cut-off subgrid terms were well-reproduced while the low-wavenumber were underestimated. For K-S, the method performed better overall. For the Burgers equation, a simple power law scaling between the size $N_c$ of the resolved set F and the memory length was observed. For K-S, the estimated memory length was shorter than for the Burgers equation. \\
\indent Further work will examine the validity of our ersatz approach. Other solution methods to solve the orthogonal dynamics equation will be investigated as well.

\section{Acknowledgments}
This research was funded by the AFOSR under the project \textit{LES Modeling of Non-local effects using Statistical Coarse-graining} (Tech.  Monitor:  Jean-Luc Cambier).

\end{document}